\documentclass[11pt]{article}
\usepackage{amsfonts}
\usepackage{amssymb,amsmath}
\usepackage{latexsym}
\usepackage{graphics}
\usepackage{epsfig}
\usepackage{psfrag}

\numberwithin{equation}{section}

\def\RR{\mathbb{R}}

\def\eps{\varepsilon}

\newtheorem{thm}{Theorem}[section]
\newtheorem{lem}[thm]{Lemma}

\newtheorem{prop}[thm]{Proposition}

\newtheorem{property}{Property}
\newenvironment{proof}{\paragraph{Proof} \phantom{9}}{\hfill$\Box$\bigskip}

\begin{document}

\title{Well-posedness in any dimension for Hamiltonian flows with
  non $BV$ force terms}

\author{Nicolas Champagnat$^{1}$, Pierre-Emmanuel Jabin$^{1,2}$}

\footnotetext[1]{TOSCA project-team, INRIA Sophia Antipolis --
  M\'editerran\'ee, 2004 rte des Lucioles, BP.\ 93, 06902 Sophia
  Antipolis Cedex, France, \\
E-mail: \texttt{Nicolas.Champagnat@sophia.inria.fr}}

\footnotetext[2]{Laboratoire J.-A. Dieudonn\'e, Universit\'e de Nice --
  Sophia Antipolis, Parc Valrose, 06108 Nice Cedex 02, France, E-mail:
  \texttt{jabin@unice.fr}}
\date{}

\maketitle

\begin{abstract}
We study existence and uniqueness for the classical dynamics of a
particle in a force field in the phase space. Through an explicit
control on the regularity of the trajectories, we show that this is
well posed if the force belongs to the Sobolev space $H^{3/4}$.
\end{abstract}
\bigskip

\bigskip


\section{Introduction}
\label{sec:intro}

This paper studies existence and uniqueness of a flow for the equation
\begin{equation}
  \label{eq:Newton}
  \begin{cases}
    \partial_t X(t,x,v)=V(t,x,v), \quad & X(0,x,v)=x, \\
    \partial_t V(t,x,v)=F(X(t,x,v)), \quad & V(0,x,v)=v,
  \end{cases}
\end{equation}
where $x$ and $v$ are in the whole $\RR^d$ and $F$ is a given
function from $\RR^d$ to $\RR^d$. Those are of course Newton's equations for a
particle moving in a force field $F$. For many applications the force
field is in fact a potential
\begin{equation}
F(x)=-\nabla\phi(x),
\end{equation} 
even though we will not use the additional Hamiltonian structure that
this is providing.

This is a particular case of a system of differential equations
\begin{equation}
\partial_t \Xi(t,\xi)=\Phi(\Xi),\label{ode}
\end{equation}
with $\Xi=(X,V)$, $\xi=(x,v)$, $\Phi(\xi)=(v,F(x))$.
 Cauchy-Lipschitz'
Theorem applies to \eqref{eq:Newton} and gives maximal solutions if
$F$ is Lipschitz. Those solutions are in particular global in
time if for instance $F\in L^\infty$. Moreover because of the
particular structure of Eq.~\eqref{eq:Newton}, this solution has the
additional 
\begin{property}
For any
  $t\in\RR$ the application
  \begin{equation}
    (x,v)\in\RR^d\times\RR^d\mapsto(X(t,x,v),V(t,x,v))\in\RR^d\times\RR^d
   \label{invertible}
  \end{equation}
  is globally invertible and has Jacobian $1$ at any
  $(x,v)\in\RR^d\times\RR^d$.
 It also defines a semi-group 
\begin{equation}\begin{aligned}
&  \forall s,t\in\RR,\qquad & X(t+s,x,v)=X(s,X(t,x,v),V(t,x,v)), \\
 & \mbox{and}\qquad & V(t+s,x,v)=V(s,X(t,x,v),V(t,x,v)).
\end{aligned}\label{semigroup}\end{equation}
\label{prop:Hamilt}\end{property}
In many cases this Lipschitz
regularity is too demanding and one would like to have a well
posedness theory with a less stringent assumption on $F$. That is the
aim of this paper. More precisely, we prove
\begin{thm}
  \label{thm:main}
  Assume that $F\in H^{3/4}\cap
  L^{\infty}$. 
Then, there exists a solution 
  to (\ref{eq:Newton}), satisfying Property \ref{prop:Hamilt}. 
Moreover this solution is unique among all limits
  of solutions to any regularization of \eqref{eq:Newton}.
\end{thm}
Many works have already studied the well posedness of Eq.~\eqref{ode}
under weak conditions for $\Phi$. The first one was essentially due to
DiPerna and Lions \cite{DL}, using the connection between \eqref{ode}
and the transport equation
\begin{equation}
\partial_t u+\Phi(\xi)\cdot\nabla_\xi u=0.\label{transport}
\end{equation}
The notion of renormalized solutions for Eq.~\eqref{transport}
provided a well posedness theory for \eqref{ode} under the conditions
$\Phi\in W^{1,1}$ and ${\rm div}_\xi \Phi\in L^\infty$. This theory
was generalized in \cite{Li2}, \cite{LL} and \cite{HLL}. 

Using a slightly different notion of renormalization, Ambrosio
\cite{Am} obtained well posedness with only $\Phi\in BV$ and
${\rm div}_\xi \Phi\in L^\infty$ (see also the papers by Colombini and
Lerner \cite{CL}, \cite{CL2} for the $BV$ case). The bounded
divergence condition was then slightly relaxed by Ambrosio, De Lellis
and Mal\`y in \cite{ADM} with then
only $\Phi\in SBV$ (see also \cite{DeL}).

Of course there is certainly a limit to how weakly $\Phi$ may be and
still provide uniqueness, as shown by the counterexamples of Aizenman
\cite{Ai} and Bressan \cite{Br1}. The example by De Pauw \cite{DeP}
even suggests that for the general setting \eqref{ode}, $BV$ is
probably close to optimal.

But as \eqref{eq:Newton} is a very special case of \eqref{ode},
it should be easier to deal with. And for instance Bouchut \cite{Bo} got
existence and uniqueness to \eqref{eq:Newton} with $F\in BV$ in a
simpler way than \cite{Am}. Hauray \cite{Ha2} handled a slightly less
than $BV$ case ($BV_{loc}$). 

In dimension $d=1$ of physical space (dimension $2$ in phase space),
Bouchut and Desvillettes proved well posedness for Hamiltonian systems
(thus including \eqref{eq:Newton} as $F$ is always a derivative in
dimension $1$) without any additional derivative for $F$ (only
continuity). This was extended to Hamiltonian systems in dimension $2$
in phase space with only $L^p$ coefficients in \cite{Ha1} and even to
any system (non necessarily Hamiltonian) 
with bounded divergence and continuous coefficient by
Colombini, Crippa and Rauch \cite{CCR} (see also \cite{CR} for low
dimensional settings and \cite{BJ} with a very different goal in
mind).

Unfortunately in large dimensions (more than $1$ of physical space or
$2$ in the phase space), the Hamiltonian or bounded divergence
structure does not help so much. To our knowledge, Th.~\ref{thm:main}
is the first result to require less than $1$ derivative on the force
field $F$ in any dimension. Note that the comparison between $H^{3/4}$
and $BV$ is not clear as obviously $BV\not\subset H^{3/4}$ and
$H^{3/4}\not\subset BV$. Even if one considers the stronger assumption
that the force field be in $L^\infty\cap BV$, that space contains by
interpolation $H^s$ for $s<1/2$ and not $H^{3/4}$. 
As the proof of Th.~\ref{thm:main} uses
orthogonality arguments, we do not know how to work in spaces non
based on $L^2$ norms ($W^{3/4,1}$ for example). Therefore strictly
speaking Th.~\ref{thm:main} is neither stronger nor weaker than
previous results.

We have no idea whether this $H^{3/4}$ is optimal or in which
sense. It is striking because it already appears in a question
concerning the related Vlasov equation
\begin{equation}
\partial_t f+v\cdot\nabla_x f+F\cdot\nabla_v f=0.\label{vlasov}
\end{equation}
Note that this is the transport equation corresponding to
Eq.~\eqref{eq:Newton}, just as Eq.~\eqref{transport} corresponds to
\eqref{ode}. As a kinetic equation, it has some regularization
property namely that the average
\[
\rho(t,x)=\int_{\RR^d} f(t,x,v)\,\psi(v)\,dv,\quad \mbox{with}\ \psi\in
C^\infty_c(\RR^d), 
\]
is more regular than $f$. And precisely if $f\in L^2$ and $F\in
L^\infty$ then $\rho\in H^{3/4}$; we refer to Golse, Lions, Perthame
and Sentis \cite{GLPS} for this result, DiPerna, Lions, Meyer for a
more general one \cite{DLM} or \cite{Ja2} for a survey of averaging
lemmas. Of course we do not know how to use this kind of result for
the uniqueness of \eqref{vlasov} or even what is the connection
between the $H^{3/4}$ of averaging lemmas and the one found here. It
{\em could} just be a scaling property of those equations.

Note in addition that the method chosen for the proof may in fact be
itself a limitation. Indeed it relies on an explicit control on the
trajectories~: for instance, we show that $|X(t,x,v)-X^\delta(t,x,v)|$
and $|V(t,x,v)-V^\delta(t,x,v)|$ remain approximately of order
$|\delta|$ if
% they were so at the initial time, for instance
%
%\begin{multline*}
%  Q_\delta(T)=\iint_{\Omega}\log\left(1+\frac{1}{|\delta|^2}\left(\sup_{0\leq
%        t\leq T}|X(t,x,v)-X^\delta(t,x,v)|^2
%      \right.\right. \\
%      \left.\left. +\int_0^T|V(t,x,v)-V^\delta(t,x,v)|^2 \:dt\right)\right)
%      \:dx \:dv, 
%\end{multline*}
\[
X^\delta(t,x,v)=X(t,x+\delta_1,v+\delta_2),\quad 
V^\delta(t,x,v)=V(t,x+\delta_1,v+\delta_2).
\]  
However the example given in Section \ref{example} demonstrates that
such a control in not always possible: Even in $1d$ it requires at
least $1/2$ derivative on the force term ($F\in W^{1/2,1}_{loc}$)
whereas well posedness is known with essentially $F\in L^p$ (see the
references above).

This kind of control is obviously connected with
regularity properties of the flow (differentiability for instance),
which were studied in \cite{ALM} (see also \cite{AC}). The idea to
prove them directly and then use them for well posedness is quite
recent, first by Crippa and De Lellis in \cite{CD} with the
introduction and subsequent bound on the functional
\begin{equation}
\int_\Omega \sup_r \int_{|\delta|\leq r}
\log\left(1+\frac{|\Xi(t,\xi)-\Xi(t,\xi+\delta)|}{|\delta|}\right)
\,d\delta\,dx. \label{Q}
\end{equation}
This gave existence/uniqueness for Eq.~\eqref{ode} with $\Phi\in
W^{1,p}_{loc}$ for any $p>1$ and a weaker
version of the bounded divergence condition. This was extended in
\cite{BC} and \cite{Ja}.

We use here a modified version of \eqref{Q} which takes the different
roles of $x$ and $v$ into account. The way of bounding it is also
quite different as we essentially try to integrate the oscillations of
$F$ along a trajectory. 

\medskip

The paper is organized as follows: The next section introduces the
functional that is studied, states the bounds that are to be proved
and briefly explains the relation with the well posedness result
Th.~\ref{thm:main}. The section after that presents the example in
$1d$ and the last and longer section the proof of the bound.

\paragraph{Notation}

\begin{itemize}
\item $u\cdot v$ denotes the usual scalar product of $u\in\RR^d$ and
  $v\in\RR^d$.
\item $S^{d-1}$ denotes the $d-1$-dimensional unit sphere in $\RR^d$.
\item $B(x,r)$ is the closed ball of $\RR^{d}$ for the standard
  Euclidean norm with center $x\in\RR^d$ and radius $r\geq 0$.
\item $C$ denotes a positive constant that may change from line to
  line.
\end{itemize}

%%%%%%%%%%%%%%%%%%%%%%%%%%%%%%%%%
\section{Preliminary results}
\label{sec:prelim}
%%%%%%%%%%%%%%%%%%%%%%%%%%%%%%%%%
\subsection{Reduction of the problem}
\label{sec:reduc}
%%%%%%%%%%%%%%%
%
In the sequel, we give estimates on the flow to Eq. (\ref{eq:Newton}) for
initial values $(x,v)$ in a compact subset
$\Omega=\Omega_1\times\Omega_2$ of $\RR^{2d}$ and for time
$t\in[0,T]$. Fix some $A>0$ and  consider any $F\in L^\infty$ with
$\|F\|_{L^\infty}\leq A$. Then for any solution to Eq.~\eqref{eq:Newton}  
\begin{equation*}\begin{split}
 |V(t,x,v)-v|  & \leq \|F\|_{L^\infty} t\leq A\,t\\
 \mbox{and}\quad
  |X(t,x,v)-x| & \leq vt +\|F\|_{L^\infty} t^2/2\leq vt+A\,t^2/2.
\end{split}
\end{equation*}
Therefore, for any $t\in[0,T]$ and for any $(x,v)$ at a distance
smaller than 1 from $\Omega$,
$(X(t,x,v),V(t,x,v))\in\Omega'=\Omega'_1\times\Omega'_2$ for some
compact subset $\Omega'$ of $\RR^{2d}$. Moreover $\Omega'$ depends
only on $\Omega$ and $A$.
Similarly, we introduce $\Omega''$ a compact subset of $\RR^{2d}$ such
that the couple $(X(-t,x,v),V(-t,x,v))$ belongs to $\Omega''$ for any
$t\in[0,T]$ and any 
$(x,v)$ at a distance smaller than 1 from $\Omega'$.

% Thanks to those definitions, we could assume that $F$ is in fact
% compactly supported in $\Omega''$.

\medskip

For $T>0$, define the quantity
\begin{multline*}
  Q_\delta(T)=\iint_{\Omega}\log\left(1+\frac{1}{|\delta|^2}\left(\sup_{0\leq
        t\leq T}|X(t,x,v)-X^\delta(t,x,v)|^2
      \right.\right. \\
      \left.\left. +\int_0^T|V(t,x,v)-V^\delta(t,x,v)|^2 \:dt\right)\right)
      \:dx \:dv, 
\end{multline*}
where $X,V$ and $X^\delta, V^\delta$ are two solutions to
\eqref{eq:Newton}, satisfying
\begin{equation}\begin{split}
& X(0,x,v)=x,\qquad V(0,x,v)=v,\\
&|X(0,x,v)-X^\delta(0,x,v)|\leq |\delta|,\qquad 
|V(0,x,v)-V^\delta(0,x,v)|\leq |\delta|.\label{initial}
\end{split}\end{equation}

We prove the following result
\begin{prop}
  \label{prop:control-Q}
  Fix $T>0$, any $A>0$ and $\Omega\in\RR^{2d}$ compact. Define
  $\Omega'$ and $\Omega''$ as in Section~\ref{sec:reduc}. There exists
  a constant $C>0$ depending only of $\mbox{diam}(\Omega')$,
  $|\Omega''|$, $T$ and $A$, such that, for any $a\in(0,1/4)$, $F\in
  H^{3/4+a}$ with $\|F\|_{L^{\infty}}\leq A$ and any solutions $(X,V)$
  and $(X^\delta,V^\delta)$ to \eqref{eq:Newton} satisfying
  Property~\ref{prop:Hamilt} and \eqref{initial}, one has for any
  $|\delta|<1/e$,
  \begin{equation*}
    Q_\delta(T)\leq
    C\left(1+\left(\log\frac{1}{|\delta|}\right)^{\max\{1-2a,1/2\}}\right)
    \left(1+\|F\|_{H^{3/4+a}(\Omega'')}\right).
  \end{equation*}
\end{prop}
 As will appear in the proof, this result can be actually
  extended without difficulty to any $F\in L^{\infty}$ such that
  $$
  \int_{\RR^d}|k|^{3/2}|\alpha(k)|^2f(k)\:dk<\infty
  $$
  for some function $f\geq 1$ such that $f(k)\rightarrow+\infty$ when
  $|k|\rightarrow+\infty$, where $\alpha(k)$ is the Fourier
  transform of $F$. We restrict ourselves to Prop.~\ref{prop:control-Q}
  to simplify the presentation in the proof
 but this remark means that the following
 modified proposition holds
\begin{prop}
  \label{prop:control-Qbis}
  Fix $T>0$, $A>0$, $\Omega\in\RR^{2d}$ compact and any $f\geq 1$ such that
  $f(k)\rightarrow+\infty$ when 
  $|k|\rightarrow+\infty$. Define 
$\Omega'$ and $\Omega''$ as in
  Section~\ref{sec:reduc}.  
There exists a continuous, decreasing function $\eps(\delta)$ with
$\eps(0)=0$
s.t. for any $F\in H^{3/4}\cap L^\infty$  with
  $\|F\|_{L^{\infty}}\leq A$, 
 for any solutions $(X,V)$ and
  $(X^\delta,V^\delta)$ to \eqref{eq:Newton} satisfying
  Property~\ref{prop:Hamilt} and \eqref{initial}, one has
  for any 
$|\delta|<1/e$,
  \begin{equation*}
    Q_\delta(T)\leq
   |\log \delta|\;\eps(\delta)\;\left(1+\int_{\RR^d} |k|^{3/2}\,
   |\alpha(k)|^2\,f(k)\,dk\right)^{1/2},
  \end{equation*}
with $\alpha$ the Fourier transform of $F$.
\end{prop}

%%%%%%%%%%%%%%%%%%%%%
\subsection{From Prop.~\ref{prop:control-Q} or \ref{prop:control-Qbis} 
to Th.~\ref{thm:main}}
%%%%%%%%%%%%%%%%%
It is now well known how to pass from an estimate like the one provided
by Prop.~\ref{prop:control-Q} to a well posedness theory (see \cite{CD}
for example) and therefore
we only briefly recall the main steps. Take any $F\in H^{3/4+a}\cap
L^\infty$.

We start by the existence of a solution. For that define $F_n$ a
regularizing sequence of $F$. Denote $X_n$, $V_n$ the solution to
\eqref{eq:Newton} with $F_n$ instead of $F$ and
$(X_n,V_n)(t=0)=(x,v)$. For any $\delta=(\delta_1,\;\delta_2)$ in
$\RR^{2d}$, put
\[
(X_n^\delta,\ V_n^\delta)(t,x,v)=(X_n,\ V_n)(t,x+\delta_1,v+\delta_2).
\]
The function $F_n$ and the solutions $(X,V)$, $(X^\delta,V^\delta)$
satisfy to all the assumptions of Prop.~\ref{prop:control-Q}, as
$F_n\in W^{1,\infty}$, using Property~\ref{prop:Hamilt}. Since $F\in
L^{\infty}\cap H^{3/4+a}$, one may choose $F_n$ uniformly
bounded in this space. The proposition then shows that
\begin{multline*}
  Q_{\delta,n}(T)=\iint_{\Omega}\log\left(1+\frac{1}{|\delta|^2}\left(\sup_{0\leq
        t\leq T}|X_n(t,x,v)-X^\delta_n(t,x,v)|^2
      \right.\right. \\
      \left.\left. +\int_0^T|V_n(t,x,v)-V^\delta_n(t,x,v)|^2 \:dt\right)\right)
      \:dx \:dv, 
\end{multline*}  
is uniformly bounded in $n$ and $\delta$ by
\[
C\left(1+\left(\log\frac{1}{|\delta|}\right)^{\max\{1-2a,1/2\}}\right).
\]
By Rellich criterion, this proves that the sequence $(X_n,V_n)$ is
compact in $L^1_{loc}(\RR_+\times\RR^{2d})$. Denote by $(X,V)$ an
extracted limit, one directly checks that $(X,V)$ is a solution to
\eqref{eq:Newton} by compactness and satisfies \eqref{invertible},
\eqref{semigroup}. Thus existence is proved.

For uniqueness, consider another solution $(X^\delta,V^\delta)$ to
\eqref{eq:Newton}, which is the limit of solutions to a regularized
equation (such as the one given by $F_n$ or by another regularizing
sequence of $F$). Then with the same argument,
$(X^\delta,V^\delta)$ also satisfies \eqref{invertible} and
\eqref{semigroup}. Moreover
\[
X(0,x,v)-X^\delta(0,x,v)=x-x=0,\quad V(0,x,v)-V^\delta(0,x,v)=v-v=0,
\] 
so that $(X,V)$ and $(X^\delta, V^\delta)$ also verify \eqref{initial}
for any $\delta\neq 0$. Applying again Prop.~\ref{prop:control-Q} and
letting $\delta$ go to $0$, one concludes that $X=X^\delta$ and $V=V^\delta$. 

Note from this sketch that one has 
uniqueness among all solutions to \eqref{eq:Newton} satisfying
\eqref{invertible} and \eqref{semigroup} and not only
those which are limit of a regularized problem. However not all
solutions to \eqref{eq:Newton} (pointwise) necessarily satisfy those
two conditions so that the uniqueness among all solutions to
\eqref{eq:Newton} is unknown. Indeed in many cases, it is not true, as
there is a hidden selection principle in \eqref{invertible} (see the
discussion in \cite{ADM}, \cite{Cr} or \cite{DeL}).

Finally if $F\in H^{3/4}$ only, then one first applies 
the De La Vall\'ee Poussin's lemma to find a function $f$ s.t.\
$f(k)\rightarrow +\infty$ when $|k|\rightarrow+\infty$ and
\begin{equation}
\int_{\RR^d} |k|^{3/2}\,f(k)\,|\alpha(k)|^2\,dk<+\infty.
\end{equation}
One proceeds as before with a regularizing sequence $F_n$ which now
has to satisfy uniformly the previous estimate. Using
Prop.~\ref{prop:control-Qbis} instead of Prop.~\ref{prop:control-Q},
the rest of the proof is identical.
%%%5%%%%%%%%%%%%%%%%%%%%%%%%%%%%%%%%%%%%%
\section{The question of optimality : An example}
\label{example}
%%%%%%%%%%%%%%%%%%%%%%%%%%%%%%%%%%%%%%%
%%%%%%%%%%%%%%%%%%%%%%%%
It is hard to know whether the condition $F\in H^{3/4}$ is optimal and
in which sense (see the short discussion in the introduction).
Instead the purpose of this section is to give a simple example
showing that $F\in W^{1/2,1}$ is a necessary condition in order to use
the method followed in this paper; namely a quantitative estimate on
$X-X^\delta$ and $V-V^\delta$. More precisely, for any $\alpha<1/2$,
we are going to construct a sequence of force fields $(F_N)_{N\geq 1}$
uniformly bounded in $W^{\alpha,1}\cap L^\infty$ and a sequence
$(\delta_N)_{N\geq 1}$ converging to 0 such that functionals like
$Q_\delta(T)$ cannot be uniformly bounded in $N$.

This example is one dimensional ($2$ in phase space) where it is known
that much less is required to have uniqueness of the flow (almost $F$
a measure). So this indicates in a sense that the method itself is
surely not optimal. Moreover what this should imply in higher
dimensions is not clear...

Through all this section we use the notation $f=O(g)$ if there exists
a constant $C$ s.t.
\[
|f|\leq C\,|g|\ \mbox{a.e.}
\]
In dimension $1$ all $F$ derive from a potential so take 
\[
\phi(x)=x+\frac{h(N\,x)}{N^{\alpha+1}},\ F=-\phi'(x)
\]
with $h$ a periodic and regular function ($C^2$ at least) with
$h(0)=0$.

As $\phi$ is regular, we know that the solution $(X,V)$ with initial
condition $(x,v)$ and the shifted one $(X^\delta,V^\delta)$
corresponding to the initial condition $(x,v+\delta)$ satisfy the
conservation of energy or
\[
V^2+2\phi(X)=v^2+2\phi(x),\quad |V^\delta|^2+2\phi(X^\delta)=|v+\delta|^2+2\phi(x).
\]  
As $\phi$ is
defined up to a constant, we do not need
to look at all the trajectories and may instead restrict ourselves to
the one starting at $x$ s.t. $v^2+\phi(x)=0$. 
By symmetry, we may assume $v>0$ and
excluding the negligible set of initial data with $v=0$, we may even
take $v>\delta$.

Let $t_0$ and $t_0^\delta$ be the first times when the trajectories
stop increasing: $V(t_0)=0$ and $V^\delta(t_0^\delta)=0$. As both
velocities are initially positive, they stay so until $t_0$ or
$t_0^\delta$. So for instance
\[
\dot X=V=\sqrt{-2\phi(X)}.
\]
Hence $t_0$ is obtained by
\[\begin{split}
t_0&=\int_0^{t_0} \frac{\dot X}{\sqrt{-2\phi(X)}}\,dt=\int_x^{x_0}
\frac{dy}{\sqrt{-2\phi(y)}}\\
&=\int_x^{x_0}
\frac{dy}{\sqrt{-2y-2h(Ny)\,N^{-1-\alpha}}},
\end{split}\]
if $x_0=X(t_0)$. Of course by energy conservation $\phi(x_0)=0$ and
again as we are in dimension $1$ this means that we may simply take $x_0=0$.

\begin{figure}
  \label{fig:ex}
  \begin{center}
    \psfrag{phi}{$\phi$} \psfrag{x=0}{$x_0=0$}
    \psfrag{xdel}{$x^\delta_0$} \psfrag{x}{$x$}
    \psfrag{v2}{$v^2$} \psfrag{v+del}{$(v+\delta)^2$}
    \includegraphics[width=10cm]{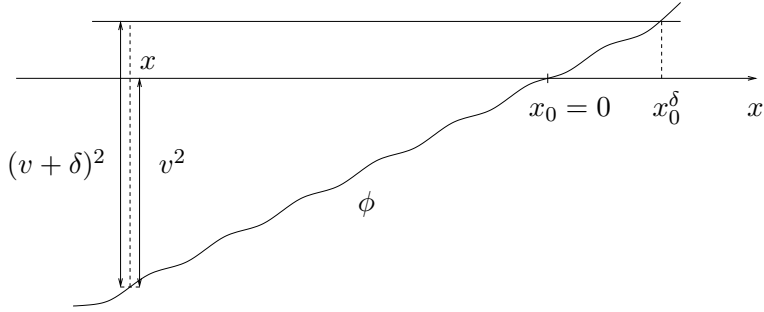}
  \end{center}
  \caption{The potential $\phi$ and the construction of $x_0$ and $x^\delta_0$}
\end{figure}

We have the equivalent formula for $t_0^\delta$ with $x_0^\delta$
(which we may not assume equal to $0$). Put
\[
C^\delta=|v+\delta|^2+2\phi(x)=\delta^2+2v\delta,\quad 
\eta=N(x_0^\delta-x_0)= N\,x_0^\delta
\]
and note that $2\phi(x_0^\delta)-2\phi(x_0)=C^\delta$, so that
$|x_0^\delta-x_0|=|x_0^\delta|\leq C\delta$ since $\phi'\geq 1/2$ for
$N$ large enough.  Then
\[\begin{split}
t_0^\delta&=\int_x^{x_0^\delta} \frac{dy}{\sqrt{C^\delta-2\phi(y)}}=
\int_{x-\eta/N}^{x_0}
\frac{dy}{\sqrt{C^\delta-2y
-2\eta/N-2\,N^{-1-\alpha}\,h(N\,y+\eta)}}\\ 
&=O(\delta)+\int_{x}^{x_0}
\frac{dy}{\sqrt{-2y-2\,N^{-1-\alpha}\,(h(N\,y+\eta)-h(\eta))}},
\end{split}\] as the integral between $x$ and $x-\eta/N$ is bounded by
$O(\delta)$ (the integrand is bounded here) and
\[
C^\delta=2\phi(x_0^\delta)=2x_0^\delta+\frac{2}{N^{1+\alpha}}h(N\,x_0^\delta)
=2\frac{\eta}{N}+\frac{2}{N^{1+\alpha}} h(\eta).
\]
Note that as $h$ is Lipschitz regular
\[\begin{split}
&\frac{|h(Nx+\eta)-h(\eta)|}{N^{1+\alpha}}=
O\left(\frac{x}{N^\alpha}\right), \\ 
&\frac{|h(Nx)|}{N^{1+\alpha}}=\frac{|h(Nx)-h(0)|}{N^{1+\alpha}}=
O\left(\frac{x}{N^\alpha}\right).
\end{split}
\]
So subtracting the two formula and making an asymptotic expansion
\[\begin{split}
t_0-t_0^\delta=O(\delta)+\int_x^{x_0}\frac{dy}{\sqrt{-2y}^3}&\left(
-\frac{2}{N^{1+\alpha}}
(h(Ny)-h(Ny+\eta)+h(\eta))\right.\\
&\left.+O\left(\frac{h(Ny)}{N^{1+\alpha}
  \sqrt{y}}\right)^2\right).
\end{split}\] Making the change of variable $Ny=z$ in the dominant
term in the integral, one finds
\[\begin{split}
t_0-t_0^\delta=&O(\delta)-2\int_{Nx}^0 N^{1/2-1-\alpha}\,
\frac{h(z)-h(z+\eta)+h(\eta)}{\sqrt{-2z}^3}\,dz \\
&+O(N^{-3/2-2\alpha}).
\end{split}\]
Consequently as long as
\[
A(\eta)=\int_{-\infty}^0 \frac{h(z)-h(z+\eta)-h(\eta)}{\sqrt{-2z}^3}\;dz
\]
is of order $1$ then $t_0-t_0^\delta$ is of order
$N^{-1/2-\alpha}$. Note that $A(\eta)$ is small when $\eta$ is, but it
is always possible to find functions $h$ s.t.\ $A(\eta)$ is of order
$1$ at least for some $\eta$. One way to see this is by observing that
% Indeed if not $A$ would have to vanish
% everywhere but for instance
\[
A'(\eta)=-\int_{-\infty}^0
\frac{h'(y+\eta)+h'(\eta)}{\sqrt{-2y}^3}\,dy
\]
cannot vanish for all $\eta$ and functions $h$. Taking $h$ such that
$A'(\eta)\geq 1$ for $\eta$ in some non-trivial interval, we can
assume that $A$ is of order $1$ for
$\eta\in[\underline{\eta},\bar{\eta}]$ for some
$\underline{\eta}<\bar{\eta}$.

Coming back to the definition of $\eta$ and $x_0^\delta$,
$\eta\in[\underline{\eta},\bar{\eta}]$ is equivalent to
\begin{equation}
\label{reldeltaN}
\delta^2+2v\delta\in \phi([\underline{\eta}/N,\bar{\eta}/N]).
\end{equation}
Using the formula for $\phi$ and the fact that $\underline{\eta}$ and
$\bar{\eta}$ are independent of $N$ or $\delta$, we find
\begin{equation*}
\delta^2+2v\delta+O(N^{-1-\alpha})\in [\underline{\eta}/N,\bar{\eta}/N].
\end{equation*}
So let us finally choose $\delta=1/N$ and denote by ${\cal V}$ the
space of initial velocities $v$ s.t.~\eqref{reldeltaN} is satisfied
for $N$ large enough. In view of the previous computation, there
exists $N_0\geq 1$ and $\gamma>0$ such that for all $v\in{\cal V}$ and
all $N\geq N_0$, 
$$
\gamma N^{-1/2-\alpha}\leq |t_0-t_0^\delta|=|t_0(v)-t_0^\delta(v)|\leq
\gamma^{-1} N^{-1/2-\alpha}.
$$

We consider now the rest of the trajectories after times $t_0$ and
$t_0^\delta$. To this aim, we denote
by $Y(t,y)$ and $ W(t,y)$ the solution of \eqref{eq:Newton} with initial
data $(y,0)$. By uniqueness
\[
X(t)=Y(t-t_0,x_0)\quad \forall t\geq t_0\quad\mbox{and}\quad
X^\delta(t)=Y(t-t_0^\delta,x_0^\delta)\quad \forall t\geq t_0^\delta.
\] 
%Assume for
%instance that $\eta_0<0$. Then ({\bf a justifier})
%\[
%Y(t,-T/N)\leq Y(t,x_0^\delta)\leq Y(t,0),
%\] 
%where $T$ is the period of $h$. On the other hand, by periodicity of
%$h$
%\[
%Y(t,-T/N)=Y(t,0)-T/N.
%\]
%This implies that
%\[
%X^\delta(t)=X(t+t_0^\delta-t_0)+O(1/N).
%\]
Obviously one cannot have $V(t)$ small for all times, as initially
$v\in {\cal V}$ was not small, and, as the force field $\nabla\phi$ is
bounded, $V$ is Lipschitz in time. So in conclusion for any $v\in
{\cal V}$, there exists a time interval $I\subset (t_0,+\infty)$ of
length of order $v$ where $V$ is larger than $v/2$.

Moreover $x_0^\delta\in
x_0+[\underline{\eta}/N,\bar{\eta}/N]=[\underline{\eta}/N,\bar{\eta}/N]$. Now
either there exists a time interval $J$ of order $v$ s.t.
\[
\forall t\in J,\quad |Y(t,x_0)-Y(t,x_0^\delta)|\geq \gamma N^{-1/2-\alpha}\,v/4. 
\]
or if it is not the case then on a subset $\tilde I$ of $I$ of size
$v$, 
one has
\[
|Y(t-t_0,x_0)-Y(t-t_0,x_0^\delta)|\leq \gamma N^{-1/2-\alpha}\,v/4.
\] 
Note that $t_0$ may be replaced by $t_0^\delta$ in the previous
inequality by reducing the interval $\tilde{I}$ (while keeping its
length of order 1) since $|t_0-t_0^\delta|=O(N^{-1/2-\alpha})=o(1)$.
Therefore for $t\in \tilde I$
\[\begin{split}
|X(t)-X^\delta(t)|&\geq
|X(t)-X(t+t_0-t_0^\delta)|-|Y(t-t_0^\delta,x_0)-Y(t-t_0^\delta,x_0^\delta)|\\
&\geq |t_0-t_0^\delta|\,v/2-\gamma N^{-1/2-\alpha}\,v/4\geq \gamma\,v\,N^{-1/2-\alpha}/4,
\end{split}\]
as $V$ is larger than $v/2$.

Consequently in both situations, we have two solutions,
$(Y(t,x_0),W(t,x_0))$ and $(Y(t,x_0^\delta),W(t,x_0^\delta))$ or
$(X,V)$ and $(X^\delta,V^\delta)$, distant of $1/N$ initially but
distant of order $N^{-1/2-\alpha}$ on a time interval of order
$v$. Since $h$ is periodic this provides many initial conditions with
such a property. The difficulty that the distance between $x_0$ and
$x_0^\delta$ is not fixed can be overcome since we are in
two-dimensional setting (another trajectory starting further than
$x_0^\delta$ from $x_0$ cannot approach more $(Y(t,x_0),W(t,x_0))$
than $(Y(t,x_0^\delta),W(t,x_0^\delta))$ does).  Therefore we may
control a functionals like $Q_\delta(T)/(\log(1/\delta))^{1-a}$ with
$a>0$ uniformly in $N$ only if
\[
N^{-1/2-\alpha}=O(\delta).
\]
Since $\delta=1/N$, this requires
\[
\alpha\geq 1/2,
\]
or $F=-\nabla\phi$ in at least $W^{1/2,1}$ as claimed.
%
%%%%%%%%%%%%%%%%%%%%%%%%%%%%%%%%%%%%%%%%%%%%%%%%%%%%%%%%%%%%%%
%%%%%%%%%%%%%%%%%%%%%%%%%%%%%%%%%%%%%%%%%%%%%%%%%%%%%%%%%%%%%%%
\section{Control of $Q_\delta(T)$~: Proof of Prop.~\ref{prop:control-Q}}
\label{sec:control}
%%%%%%%%%%%%%%%%%%%%%%%%%%%%%%%%%%%%%%%%%%%%%%%%%%%%%%%%%%%%%%%%
%%%%%%%%%%%%%%%%%%%%%%%%%%%%%%%%%%%%%%%%%%%%%%%%%%%%%%%%%%%%%%
Recall the notation $\alpha$ for the Fourier transform of $F$. The
assumption of Proposition~\ref{prop:control-Q} corresponds to the
following bound:
\begin{equation*}
%   \label{eq:hyp-alpha}
  \int_{\RR^d}|k|^{\frac{3}{2}+2a}|\alpha(k)|^2\:dk
  =\|F\|_{H^{3/4+a}(\Omega'')}^2<+\infty.
\end{equation*}
\subsection{Decomposition of $Q_\delta(T)$}
\label{sec:decomp1}

Let
\begin{equation*}\begin{split}
  A_\delta(t,x,v)=&|\delta|^2+\sup_{0\leq
    s\leq t}|X(s,x,v)-X^\delta(s,x,v)|^2\\
&  +\int_0^t|V(s,x,v)-V^\delta(s,x,v)|^2 \:ds.
\end{split}\end{equation*}
From~(\ref{eq:Newton}), we compute
\begin{align*}
  & \frac{d}{dt}\log\Biggl(1+\frac{1}{|\delta|^2}\Big(\sup_{0\leq s\leq
        t}|X(s,x,v)-X^\delta(s,x,v)|^2\\
&\qquad\qquad\qquad\qquad\qquad  
+\int_0^t|V(s,x,v)-V^\delta(s,x,v)|^2 \:ds\Big)\Bigg) \\
  & \qquad 
  =\frac{2}{A_\delta(t,x,v)}\Bigg(\frac{d}{dt}\left(\sup_{0\leq s\leq
      t}|X(s,x,v)-X^\delta(s,x,v)|^2\right) \\ &
%  -\frac{2}{A_\delta(t,x,v)}
(V(t,x,v)-V^\delta(t,x,v))\int_0^t(F(X(s,x,v))
  -F(X^\delta(s,x,v))) \:ds\Bigg)
\end{align*}
Since, for any $f\in BV$,
\begin{equation*}
  \frac{d}{dt}\left(\max_{0\leq s\leq s} f(s)^2\right)
  \leq 2|f(s)f'(s)|\leq 4|f(s)|^2+4|f'(s)|^2, 
\end{equation*}
we deduce from the previous computation that
\begin{align*}
  Q_\delta(T) & \leq 4\iint_\Omega\int_0^T
  \frac{|X-X^\delta|^2+|V-V^\delta|^2}
  {A_\delta(t,x,v)} \:dt \:dx \:dv + \tilde{Q}_\delta(T) \\ & \leq
  4|\Omega|(1+T)+\tilde{Q}_\delta(T)
\end{align*}
where,
\begin{multline*}
  \tilde{Q}_\delta(T)=-2\int_0^T\iint_{\Omega}
  \frac{V(t,x,v)-V^\delta(t,x,v)}{A_\delta(t,x,v)}\:\cdot \\
  \int_{0}^{t}\int_{\RR^d}\alpha(k)\left(e^{ik\cdot X(s,x,v)}-e^{ik\cdot
    X^\delta(s,x,v)}\right) \:dk \:ds \:dx \:dv \:dt.
\end{multline*}

We introduce a $C_b^{\infty}$ function $\chi:\RR_+\rightarrow[0,1]$
such that $\chi(x)=0$ if $x\leq 1$ and $\chi(x)=1$ if $x\geq 2$.
Writing $X_t$ (resp.\ $V_t$) for $X(t,x,v)$ (resp.\ $V(t,x,v)$) and
$X^\delta_t$ (resp.\ $V^\delta_t$) for $X^\delta(t,x,v)$ (resp.\
$V^\delta(t,x,v)$), and introducing
\begin{equation*}
  \tilde{\alpha}(k)=
  \begin{cases}
    \alpha(k) & \mbox{if\ }|k|\geq(\log 1/|\delta|)^{2} \\
    0 & \mbox{otherwise,}    
  \end{cases}
\end{equation*}
we may write
\begin{equation*}
  \tilde{Q}_\delta(T)=\tilde{Q}^{(1)}_\delta(T)+\tilde{Q}^{(2)}_\delta(T)
  +\tilde{Q}^{(3)}_\delta(T)+\tilde{Q}^{(4)}(T),
\end{equation*}
where
\begin{multline*}
  \tilde{Q}^{(1)}_\delta(T)=-2\int_0^T\iint_{\Omega}\int_0^t
  \chi\left(\frac{|X_s-X^\delta_s|}{|\delta|^{4/3}}\right)
  \frac{V_t-V_t^\delta}{A_\delta(t,x,v)}\:\cdot \\
  \int_{\RR^d}\tilde{\alpha}(k)\left(e^{ik\cdot
      X_s}-e^{ik\cdot X^\delta_s}\right) \:dk \:ds \:dx
  \:dv \:dt,
\end{multline*}
\begin{multline*}
  \tilde{Q}^{(2)}_\delta(T)=-2\int_0^T\iint_{\Omega}\int_0^t
  \chi\left(\frac{|X_s-X^\delta_s|}{|\delta|^{4/3}}\right)
  \frac{V_t-V_t^\delta}{A_\delta(t,x,v)}\:\cdot \\ \int_{\RR^d}
  (\alpha(k)-\tilde{\alpha}(k))
\left(e^{ik\cdot X_s}-e^{ik\cdot X^\delta_s}\right) \:dk
  \:ds \:dx \:dv \:dt,
\end{multline*}
\begin{multline*}
  \tilde{Q}^{(3)}_\delta(T)=-2\int_0^T\iint_{\Omega}\int_0^{t}
  \left(1-\chi\left(\frac{|X_s-X^\delta_s|}{|\delta|^{4/3}}\right)\right)
  \frac{V_t-V_t^\delta}{A_\delta(t,x,v)}\:\cdot \\
  \int_{\{|k|\leq|\delta|^{-4/3}\}}
  \alpha(k)\left(e^{ik\cdot X_s}-e^{ik\cdot X^\delta_s}\right) \:dk \:ds \:dx
  \:dv \:dt,
\end{multline*}
and
\begin{multline*}
  \tilde{Q}^{(4)}_\delta(T)=-2\int_0^T\iint_{\Omega}\int_{0}^{t}
  \left(1-\chi\left(\frac{|X_s-X^\delta_s|}{|\delta|^{4/3}}\right)\right)  
  \frac{V_t-V_t^\delta}{A_\delta(t,x,v)}\:\cdot \\
  \int_{\{|k|>|\delta|^{-4/3}\}}
  \alpha(k)\left(e^{ik\cdot X_s}-e^{ik\cdot X^\delta_s}\right) \:dk \:ds \:dx
  \:dv \:dt.
\end{multline*}

The proof is based on a control each of these terms. As proved in
Subsection~\ref{sec:Q^4}, the fourth term can be bounded with
elementary computations. In Subsection~\ref{sec:Q^2-Q^3}, the second
and third terms are bounded using standard results on maximal
functions. Finally, the control of $\tilde{Q}_\delta^{(1)}(T)$
requires a more precise version of the maximal inequality, detailed in
Subsection~\ref{sec:Q^1}.

\subsection{Control of $\tilde{Q}^{(4)}_\delta(T)$}
\label{sec:Q^4}

Let us first state and prove a result that is used repeatedly in the
sequel.
\begin{lem}
  \label{lem:1}
  There exists a constant $C$ such that, for $|\delta|$ small enough,
  \begin{equation*}
    \int_s^T \frac{|V_t-V_t^\delta|}{\sqrt{A_\delta(t,x,v)}}\:dt\leq
    C(\log 1/|\delta|)^{1/2}. 
  \end{equation*}
\end{lem}

\begin{proof}
  Using Cauchy-Schwartz inequality,
  \begin{align*}
    \int_s^T \frac{|V_t-V_t^\delta|}{\sqrt{A_\delta(t,x,v)}}\:dt
    & \leq\int_s^T\frac{|V_{t}-V^\delta_{t}|}{\left(|\delta|^2+
        \int_0^t|{V}_{r}-V^\delta_{r}|^2\:dr\right)^{1/2}}\:dt \\
    & \leq \sqrt{T}\left(
      \int_s^T\frac{|{V}_{t}-V^\delta_{t}|^{2}}{|\delta|^2+
        \int_0^t|{V}_{r}-V^\delta_{r}|^2\:dr}\:dt
    \right)^{1/2} \\ & =\sqrt{T}\left(\log\left(
        \frac{|\delta|^2+
          \int_0^T|{V}_{r}-V^\delta_{r}|^2\:dr}{|\delta|^2+
          \int_0^s|{V}_{r}-V^\delta_{r}|^2\:dr}\right)\right)^{1/2}
    \\ & \leq C\sqrt{T}\left(\log 1/|\delta|\right)^{1/2}
  \end{align*}
  for $|\delta|$ small enough.
\end{proof}

Let us define the function
$$
\tilde{F}(x)=\int_{\{|k|>|\delta|^{-4/3}\}}\alpha(k)e^{ik\cdot x}
\:dx.
$$
Since $\sqrt{A_\delta(t,x,v)}\geq|\delta|$, we have
\begin{align*}
  |\tilde{Q}^{(4)}_\delta(T)| & \leq C\int_0^T\iint_{\Omega}
  (|\tilde{F}(X_s)|+|\tilde{F}(X^\delta_s)|)\\
&\qquad\qquad\qquad \times \int_{s}^T\frac{|V_t-V_t^\delta|}
  {|\delta|\sqrt{A_\delta(t,x,v)}}\:dt\:dx\:dv\:ds. \notag \\
  & \leq C(\log 1/|\delta|)^{1/2}|\delta|^{-1}\int_0^T\iint_{\Omega'}
  |\tilde{F}(x)|\:dx\:dv\:ds \notag \\ & \leq
  C(\log 1/|\delta|)^{1/2}|\delta|^{-1}\left(\int_{\Omega_1'}
    |\tilde{F}(x)|^2\:dx\right)^{1/2}, \notag
\end{align*}
where the second line follows from Lemma~\ref{lem:1} and from
Property~\ref{prop:Hamilt} applied to the change of variables
$(x,v)=(X_s,V_s)$ and $(x,v)=(X^\delta_s,V^\delta_s)$. Then, it
follows from Plancherel's identity that
\begin{align*}
  |\tilde{Q}^{(4)}_\delta(T)| & \leq C
  (\log 1/|\delta|)^{1/2}|\delta|^{-1}\left(\int_{\{|k|>|\delta|^{-4/3}\}}
    |\alpha(k)|^2\:dk\right)^{1/2} \notag \\ & \leq C
  (\log 1/|\delta|)^{1/2}|\delta|^{4a/3}\left(\int_{\RR^d}
    |k|^{\frac{3}{2}+2a}|\alpha(k)|^2\:dk\right)^{1/2}. %\label{eq:calc}
\end{align*}
%
%%%%%%%%%%%%%%%%%%%%%%%%%%%%%%%
\subsection{Control of $\tilde{Q}^{(2)}_\delta(T)$ and
  $\tilde{Q}^{(3)}_\delta(T)$}
\label{sec:Q^2-Q^3}
%%%%%%%%%%%%%%%%%%%%%%%%%%%%%%%%
We recall that the maximal function $Mf$ of $f\in L^{p}(\RR^{d})$,
$1\leq p\leq+\infty$, is defined by
\begin{equation*}
  Mf(x)=\sup_{r>0}\frac{C_d}{r^d}\int_{B(x,r)}f(z)\:dz,\qquad \forall x\in\RR^d.
\end{equation*}
We are going to use the following classical results (see \cite{St}). 
First, there exists a constant $C$ such that, for all $x,y\in\RR^d$
and $f\in L^{p}(\RR^{d})$,
\begin{equation}
  |f(x)-f(y)|\leq C\, |x-y|(M|\nabla f|(x)+M|\nabla f|(y)).
\label{maximal}
\end{equation}
Second, for all $1<p<\infty$, the operator $M$ is a linear
continuous application from $L^p(\RR^d)$ to itself.

We begin with the control of $\tilde{Q}^{(3)}_\delta(T)$. Let
\begin{equation*}
  \hat{F}(x)=\int_{\{|k|\leq|\delta|^{-4/3}\}}\alpha(k)e^{ik\cdot
    x} \:dx.
\end{equation*}
It follows from the previous inequality that
\begin{equation*}\begin{split}
  \Bigg|\int_{\{|k|\leq |\delta|^{-4/3}\}}\alpha(k)&(e^{ik\cdot
      X_s}-e^{ik\cdot X^\delta_s})\:dk\Bigg|
=|\hat F(X_s)-\hat F(X_s^\delta)|\\
&  \leq |X_s-X^\delta_s|\big(M|\nabla\hat{F}|(X_s)
  +M|\nabla\hat{F}|(X^\delta_s)\big).
\end{split}\end{equation*}
Therefore, since $1-\chi(x)=0$ if $|x|\geq 2$, following 
the same steps as for the contro, of $\tilde{Q}^{(4)}_\delta(T)$,
\begin{align*}
  |\tilde{Q}^{(3)}_\delta(T)| & \leq C\int_0^T\iint_{\Omega}\int_s^T
  \frac{|V_t-V_t^\delta|}{|\delta|\sqrt{A_\delta(t,x,v)}}\:|\delta|^{4/3}
  \\
&\qquad\qquad\qquad\qquad\big(M|\nabla\hat{F}|(X_s)+M|\nabla\hat{F}|(X^\delta_s)\big)
  \:dt\:dx\:dv\:ds.  \\ & \leq C(\log 1/|\delta|)^{1/2}|\delta|^{1/3}
  \left(\int_{\Omega'_1}(M|\nabla\hat{F}|(x))^2\:dx\right)^{1/2}
  \\ & \leq C(\log 1/|\delta|)^{1/2}|\delta|^{1/3}
  \left(\int_{\Omega'_1}|\nabla\hat{F}|^2(x)\right)^{1/2}.
\end{align*}
Then
\begin{align*}
  |\tilde{Q}^{(3)}_\delta(T)| & 
  \leq C(\log 1/|\delta|)^{1/2}|\delta|^{1/3}
  \left(\int_{\{|k|\leq|\delta|^{-4/3}\}}|k|^2|\alpha(k)|^2\:dk\right)^{1/2}
  \\ & \leq C(\log 1/|\delta|)^{1/2}|\delta|^{4a/3}
  \left(\int_{\RR^d}|k|^{\frac{3}{2}+2a}|\alpha(k)|^2\:dk\right)^{1/2}.
\end{align*}

The control of $\tilde{Q}^{(2)}_\delta(T)$ follows from a similar
computation: introducing
$F_0(x)=\int_{\{k<(\log 1/|\delta|)^{2}\}}\alpha(k)e^{ik\cdot
  x}\:dx$, we obtain
\begin{multline*}
  |\tilde{Q}^{(2)}_\delta(T)|\leq C\int_0^T\iint_{\Omega}\int_s^T
  \frac{|V_t-V_t^\delta|}{\sqrt{A_\delta(t,x,v)}}
  \:\frac{|X_s-X_s^\delta|}{\sqrt{A_\delta(t,x,v)}}
  \\
\big(M|\nabla F_0|(X_s)+M|\nabla F_0|(X^\delta_s)\big)\:dt\:dx\:dv\:dt.
\end{multline*}
Since $|X_s-X^\delta_s|\leq\sqrt{A_\delta(t,x,v)}$ for all $s\leq t$
\begin{align*}
  |\tilde{Q}^{(2)}_\delta(T)|
%   \\ & \leq C(\log 1/|\delta|)^{1/2}\int_0^T
%   \left(\iint_\Omega\big((M|\nabla^2\phi_0|(X_s))^2
%     +(M|\nabla^2\phi_0|(X^\delta))^2\big)\:dx\:dv\right)^{1/2}\:ds 
   & \leq C(\log 1/|\delta|)^{1/2}\int_0^T
  \left(\iint_{\Omega'}\big(M|\nabla F_0|(x)\big)^2\:dx\:dv\right)^{1/2}\:ds
%   \\ & \leq C(\log 1/|\delta|)^{1/2}
%   \left(\int_{\Omega'_1}|\nabla^2\phi_0|^2(x)\:dx\right)^{1/2} 
  \\ & \leq C(\log 1/|\delta|)^{1/2}
  \left(\int_{\{|k|<(\log 1/|\delta|)^{2}\}}|k|^2|\alpha(k)|^2\:dk\right)^{1/2}
  \\ & \leq C(\log 1/|\delta|)^{1-2a}
  \left(\int_{\RR^d}|k|^{\frac{3}{2}+2a}|\alpha(k)|^2\:dk\right)^{1/2}.
%   \left(\int_{\RR^d}|k|^{\frac{7}{2}+2a}|\alpha(k)|^2\:dk\right)^{1/2}\leq
%   C(\log 1/|\delta|)^{3/4},
\end{align*}
%%%%%%%%%%%%%%%%%%%%%%%%%%
\subsection{Control of $\tilde{Q}^{(1)}_\delta(T)$}
\label{sec:Q^1}
%%%%%%%%%%%%%%%%%%%%%%%%%%
The inequality \eqref{maximal} is insufficient to control
$\tilde{Q}^{(1)}_\delta(T)$. Our estimate relies on a more precise
version of this inequality, detailed below.

\subsubsection{Definition of $X^{\theta,h}_s$}
\label{sec:def-X^theta}

\begin{figure}
  \label{fig:def_Xtheta}
  \begin{center}
    \psfrag{X}{$X$} \psfrag{Xdel}{$X^\delta$}
    \psfrag{h}{$h$} \psfrag{setX}{$\{X^{\theta,h},0\leq\theta\leq 1\}$}
    \includegraphics[width=4cm]{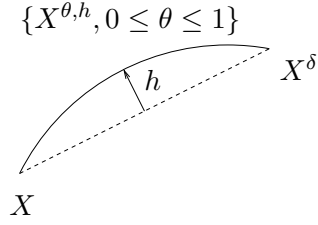}
  \end{center}
  \caption{The graph of $\theta\mapsto X^{\theta,h}$}
\end{figure}

For any $\theta\in[0,1]$ and $h\in\RR^d$, we define
\begin{equation*}
  X^{\theta,h}(t,x,v)=\theta
  X(t,x,v)+(1-\theta)X^\delta(t,x,v)
  +(1-(2\theta-1)^2)h,
\end{equation*}
and we write for simplicity $X^{\theta,h}_t$ for
$X^{\theta,h}(t,x,v)$.  Then, for any fixed $h\in\RR^d$, by
differentiation in $\theta$ 
\begin{multline}
  \int_{\RR^d}\tilde{\alpha}(k)\left(e^{ik\cdot X_s}-e^{ik\cdot
      X^\delta_s}\right) \:dk \\
  =\int_{\RR^d}\tilde{\alpha}(k)\int_0^{1}e^{ik\cdot 
    X^{\theta,h}_s}k\cdot(X_s-X^{\delta}_s+4(1-2\theta)h)
  \:d\theta \:dk.
  \label{eq:theta}
\end{multline}

For any $x,y\in\RR^d$, we introduce the hyperplane orthogonal to $x-y$ 
\begin{equation*}
  H(x,y)=\{h\in\RR^d:h\cdot(x-y)=0\}.
\end{equation*}
If $x=y$, we define for example $H(x,y)=H(0,e_1)$, where
$e_1=(1,0,\ldots,0)$.  Fix a $C^\infty_b$ function
$\psi:\RR_+\rightarrow\RR_+$ such that $\psi(x)=0$ for
$x\not\in[-1,1]$ and $\int_{H(0,e_1)}\psi(|h|) \:dh=1$. By invariance
of $|h|$ with respect to rotations, we also have
\begin{equation*}
  \int_{H(x,y)}\psi(|h|) \:dh=1
\end{equation*}
for all $x,y\in\RR^d$.

Since the left-hand side of~(\ref{eq:theta}) does not depend on $h$,
we have
\begin{multline*}
  \int_{\RR^d}\tilde{\alpha}(k)\left(e^{ik\cdot X_s}-e^{ik\cdot
      X^\delta_s}\right) \:dk \\ 
  \begin{aligned}
    =\frac{1}{|X-X^\delta|^{d-1}}
    \int_{H(X_s,X^\delta_s)}\psi\left(\frac{|h|}{|X-X^\delta|}\right)
    \int_{\RR^d}\tilde{\alpha}(k) \\ \int_0^{1}e^{ik\cdot
      X^{\theta,h}_s}k\cdot(X_s-X_s^{\delta}+4(1-2\theta)h) \:d\theta \:dk \:dh
  \end{aligned}
\end{multline*}
in the case where $X_s\not=X^\delta_s$. If $X_s=X^\delta_s$,
the previous quantity is $0$.

Let $\rho:[0,1]\rightarrow\RR_+$ be a $C^\infty_b$ function such that
$\rho(x)=1$ for $0\leq x\leq 1/4$, $\rho(x)=0$ for $3/4\leq x\leq 1$ and
$\rho(x)+\rho(1-x)=1$ for $0\leq x\leq 1$. Then, one has
\begin{equation*}
  \int_{\RR^d}\tilde{\alpha}(k)\left(e^{ik\cdot X_s}-e^{ik\cdot
      X^\delta_s}\right) \:dk= B_\delta(s,x,v)+C_\delta(s,x,v),
\end{equation*}
where
\begin{multline}
  \label{eq:def-B}
  B_\delta(s,x,v)=\frac{1}{|X_s-X_s^\delta|^{d-1}}
  \int_{H(X_s,X_s^\delta)}\psi\left(\frac{|h|}{|X_s-X_s^\delta|}\right)
  \int_{\RR^d}\tilde{\alpha}(k) \\ \int_0^{1}\rho(\theta)e^{ik\cdot
    X^{\theta,h}_s}k\cdot(X_s-X_s^{\delta}+4(1-2\theta)h) \:d\theta \:dk \:dh
\end{multline}
and
\begin{multline}
  \label{eq:def-C}
  C_\delta(s,x,v)=\frac{1}{|X_s-X_s^\delta|^{d-1}}
  \int_{H(X_s,X_s^\delta)}\psi\left(\frac{|h|}{|X_s-X_s^\delta|}\right)
  \int_{\RR^d}\tilde{\alpha}(k) \\ \int_0^{1}\rho(1-\theta)e^{ik\cdot
    X^{\theta,h}_s}k\cdot(X_s-X_s^{\delta}+4(1-2\theta)h) \:d\theta \:dk \:dh.
\end{multline}

We focus on $B_\delta(s,x,v)$ as by symmetry between $X$ and
$X^\delta$, $C_\delta$ is dealt with in exactly the same manner.
%%%%%%%%%%%%%%%%%%%%%%%%%%%%%
\subsubsection{Change of variable $z=X^{\theta,h}_s$}
\label{sec:CV-z}
%%%%%%%%%%%%%%%%%%%%%%%%%%%%%%%%%%%
For any $x\in\RR^d$, we introduce
\begin{multline}
  \label{eq:def-cone}
  K(x)=\{y\in\RR^d:\exists\theta\in[0,1],\:h\in
  H(x,0)\ \mbox{\ s.t.\ }|h|\leq|x|\\
\mbox{\ and\ }y=\theta
  (x+4(1-\theta)h)\}.
\end{multline}
\begin{figure}
  \label{fig:def_K}
  \begin{center}
    \psfrag{0}{$0$} \psfrag{x}{$x$}
    \psfrag{Nx}{$|x|$} \psfrag{K}{$K(x)$}
    \includegraphics[width=3cm]{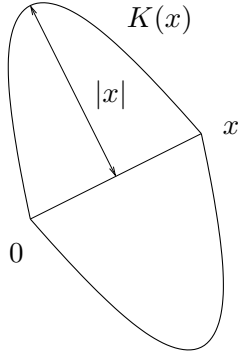}
  \end{center}
  \caption{The set $K(x)$}
\end{figure}
Observing that
\begin{equation*}
  \theta=\frac{y}{|x|}\cdot\frac{x}{|x|},
\end{equation*}
this set may also be defined as
\begin{multline*}
  K(x)=\left\{y\in\RR^d:\frac{y}{|x|}\cdot\frac{x}{|x|}\in[0,1]
    \right. \\ \left.\mbox{\ and\ } \left|\frac{y}{|x|}
      -\left(\frac{y}{|x|}\cdot\frac{x}{|x|}\right)
      \frac{x}{|x|}\right|\leq
    4\frac{y}{|x|}\cdot\frac{x}{|x|}
    \left(1-\frac{y}{|x|}\cdot\frac{x}{|x|}\right)\right\}.
\end{multline*}
Note that, for any $y\in K(x)$, taking $\theta$ and $h$ as
in~(\ref{eq:def-cone}), we have
$|y|^2=\theta^{2}(|x|^2+16(1-\theta^{2})|h|^2)\leq
17\theta^{2}|x|^2$. Therefore, denoting by $(x,y)$ the angle between
the vectors $x$ and $y$,
\begin{equation}
  \label{eq:angle}
  \cos(x,y)=\frac{x}{|x|}\cdot\frac{y}{|y|}=\frac{\theta|x|}{|y|}\geq
  17^{-1/2}.
\end{equation}

For fixed $x,y\in\RR^d$, we now introduce the application
\begin{align*}
  F_{x,y}:\  & [0,1]\times\{h\in H(x,y): |h|\leq|y-x|\} \rightarrow K(x-y) \\
  & (\theta,h)\mapsto \theta (x-y+4(1-\theta)h).
\end{align*}
It is elementary to check that $F_{x,y}$ is a bijection when $x\not=
y$, with inverse
\begin{equation*}
  F^{-1}_{x,y}(z)=\left(\frac{z}{|x-y|}\cdot\frac{x-y}{|x-y|}\:,\:
    \frac{z-\left(\frac{z}{|x-y|}\cdot\frac{x-y}{|x-y|}\right)
      (x-y)}{4\:\frac{z}{|x-y|}\cdot\frac{x-y}{|x-y|}
      \left(1-\frac{z}{|x-y|}\cdot\frac{x-y}{|x-y|}\right)}
    \right)
\end{equation*}
for $z\in K(x-y)$. Moreover, $F_{x,y}$ is differentiable and its
differential, written in a basis of $\RR^d$ with first vector
$(x-y)/|x-y|$, is
\begin{equation*}
  \nabla F_{x,y}(\theta,h)=\left(
    \begin{array}{cc}
      |x-y| & 4(1-2\theta)h \\ 0 & 4\theta(1-\theta)\mbox{Id}
    \end{array}
  \right).
\end{equation*}
Therefore, the Jacobian of $F_{x,y}$ at $(\theta,h)$ is
$(4\theta(1-\theta))^{d-1}|x-y|$.
\medskip

Making the change of variable $z=F_{X_s,X_s^\delta}(\theta,h)$ in~(\ref{eq:def-B}),
we can now compute
\begin{align}
  & B_\delta(s,x,v) \notag \\ &
  =\int_{\RR^d}\tilde{\alpha}(k)\int_0^{1}\int_{H(X_s,X_s^\delta)}
  \frac{\rho(\theta)\psi\left(\frac{|h|}{|X_s-X_s^\delta|}\right)}
  {|X_s-X_s^\delta|^d(4\theta(1-\theta))^{d-1}}
  \: e^{ik\cdot X_s^{\theta,h}} \notag \\ &\qquad\qquad
  k\cdot(X_s-X_s^\delta+4(1-\theta)h-4\theta h)
  \:(4\theta(1-\theta))^{d-1}|X_s-X_s^\delta|
  \:dh \:d\theta \:dk  \notag
%   \\ & =\int_{\RR^d}
%   k\tilde{\alpha}(k)\int_{C(X-X^\delta)}\frac{1}{|z|^d}
%   \tilde{\psi}\left(\frac{z}{|X-X^\delta|},\frac{X-X^\delta}{|X-X^\delta|}\right)
%   e^{ik\cdot(X^\delta+z)}\:k\cdot z \:dz \:dk 
  \\ & =B^1_\delta(s,x,v)-B^2_\delta(s,x,v), \label{eq:def-B1-B2}
\end{align}
with
\[\begin{split}
 B^1_\delta(s,x,v) =\int_{\RR^d}
  \tilde{\alpha}(k)\int_{\RR^d}\frac{k\cdot z}{|z|^d}
  \psi^{(1)}&\left(\frac{z}{|z|},
    \frac{X_s-X_s^\delta}{|X_s-X_s^\delta|},
    \frac{|z|}{|X_s-X_s^\delta|}\right)\\
&\qquad \qquad\qquad\qquad
  e^{ik\cdot (X_s^\delta+z)} \:dz \:dk, \notag
\end{split}
\]
and
 \[\begin{split}
B^2_\delta(s,x,v)= -\int_{\RR^d}
  \tilde{\alpha}(k)\int_{\RR^d}\frac{k}{|z|^{d-1}}
  \cdot\psi^{(2)}&\left(\frac{z}{|z|},
    \frac{X_s-X_s^\delta}{|X_s-X_s^\delta|},
    \frac{|z|}{|X_s-X_s^\delta|}\right)\\
 &\qquad\qquad\qquad e^{ik\cdot(X_s^\delta+z)}\:dz \:dk. 
\end{split}
\]
We defined, for $(a,b,c)\in S^{d-1}\times S^{d-1}\times(\RR\setminus\{0\})$,
\begin{equation*}
  \psi^{(1)}(a,b,c)=\frac{\tilde{\rho}((a\cdot
    b)c)\psi\left(\frac{|a-(a\cdot b)b|}{4(a\cdot b)(1-(a\cdot b)c)}\right)}
  {4^{d-1}(a\cdot b)^{d}(1-(a\cdot b)c)^{d-1}}
\end{equation*}
and
\begin{equation*}
  \psi^{(2)}(a,b,c)=\frac{\tilde{\rho}((a\cdot
    b)c)\psi\left(\frac{|a-(a\cdot b)b|}{4(a\cdot b)(1-(a\cdot b)c)}\right)}
  {4^{d-1}(a\cdot b)^{d-1}(1-(a\cdot b)c)^{d}}\:c(a-(a\cdot b)b),
\end{equation*}
where $\tilde{\rho}(x)=\rho(x)$ if $x\in[0,1]$, and
$\tilde{\rho}(x)=0$ otherwise.

It follows from~(\ref{eq:angle}) and from the definition of $\rho$
that these two functions have support in 
\begin{equation*}
%   \label{eq:support}
  \{(u,v)\in
  (S^{d-1})^2:\cos(u,v)\geq 17^{-1/2}\}\times [0,3/4].  
\end{equation*}
Moreover, they belong to
$C^{0,\infty,\infty}_b(S^{d-1},S^{d-1},\RR\setminus\{0\})$. Indeed,
since $\tilde{\rho}(x)=0$ for $x\geq 3/4$, the terms $(1-(a\cdot b)c)$
in the denominators do not cause any regularity problem. Moreover,
since $\psi(x)=0$ for $x\not\in[-1,1]$ and
\begin{equation*}
  \frac{|a-(a\cdot b)b|}{|a\cdot b|}\geq\frac{1}{|a\cdot b|}-1
\end{equation*}
for all $a,b\in S^{d-1}$, the terms $a\cdot b$ in the denominators do
not cause any worry either. Finally, since $\tilde{\rho}\in
C^\infty_b(\RR\setminus\{0\})$, the discontinuity of
$\tilde{\rho}$ at $0$ can only cause a problem in the neighborhood of
points such that $a\cdot b=0$ ($c$ being nonzero). Therefore, the previous
observation also solves this difficulty.

\subsubsection{Decomposition of $B^{1}_\delta(s,x,v)$: integration by
  parts}
\label{sec:IPP}

Writing $\psi^{(1)}_t$ for
\begin{equation}
  \label{eq:def-psi1_s}
  \psi^{(1)}\left(\frac{z}{|z|},\frac{X_t-X_t^\delta}{|X_t-X_t^\delta|},
    \frac{|z|}{|X_t-X_t^\delta|}\right),
\end{equation}
we decompose $B^1_\delta(s,x,v)$
\begin{align*}
  B^1_\delta(s,x,v) & =\int_{\RR^d}\int_{\RR^d}\tilde{\alpha}(k)
  \frac{e^{ik\cdot(X_s^\delta+z)}\:k\cdot z}{|z|^d}
  \psi^{(1)}_s
  \frac{i\frac{k}{|k|}\cdot V_s^\delta}
  {|k|^{-1/2}+i\frac{k}{|k|}\cdot V_s^\delta} \:dk \:dz
  \\ & +\int_{\RR^d}\int_{\RR^d}\tilde{\alpha}(k)
  \frac{e^{ik\cdot(X_s^\delta+z)}\:k\cdot z}{|z|^d}
  \psi^{(1)}_s
  \frac{|k|^{-1/2}}{|k|^{-1/2}+i\frac{k}{|k|}\cdot V_s^\delta} \:dk \:dz \\
  & =:B^{11}_\delta(s,x,v)+B^{12}_\delta(s,x,v).
\end{align*}
Now, let us write $\chi_s$ for
\begin{equation}
  \label{eq:def-chi_s}
  \chi\left(\frac{|X_s-X^\delta_s|}{|\delta|^{4/3}}\right),
\end{equation}
and let us define similarly as in~(\ref{eq:def-psi1_s})
and~(\ref{eq:def-chi_s}) the notation $\nabla_2\psi^{(1)}_s$,
$\nabla_3\psi^{(1)}_s$ and $\chi'_s$. Note that the term
$i\frac{k}{|k|}\cdot V_s^\delta\; e^{ik\cdot(X_s^\delta+z)}$ is
exactly the time derivative of $\frac{1}{|k|}\,
e^{ik\cdot(X_s^\delta+z)}$. So integrating by parts in
time, we obtain
\[
   \int_0^t \chi_s B_\delta^{11}(s,x,v)ds  
= \mbox{I}(t,x,v)-\mbox{II}(t,x,v)-\mbox{III}(t,x,v)
  -\mbox{IV}(t,x,v)-\mbox{V}(t,x,v),
\]
with
\[
\mbox{I}(t,x,v) =
  \int_{\RR^d}\int_{\RR^d}\tilde{\alpha}(k)\:\frac{k\cdot
    z}{|k|\:|z|^{d}}\:\frac{e^{ik\cdot
      (X^\delta_t+z)}\chi_t\psi^{(1)}_t}
  {|k|^{-1/2}+i\frac{k}{|k|}\cdot V^\delta_t}
  \:dk\:dz,
\]
\[
\mbox{II}(t,x,v)= \int_{\RR^d}\int_{\RR^d}\tilde{\alpha}(k)\:\frac{k\cdot
    z}{|k|\:|z|^{d}}\:\frac{e^{ik\cdot
      (x+\delta_1+z)}\chi_0\psi^{(1)}_0}
  {|k|^{-1/2}+i\frac{k}{|k|}\cdot(v+\delta_2)}
  \:dk\:dz,
\]
\[\begin{split}
\mbox{III}(t,x,v)=   \int_0^t\int_{\RR^d}\int_{\RR^d}\tilde{\alpha}(k)\:\frac{k\cdot
    z}{|k|\:|z|^{d}}\:&\frac{e^{ik\cdot(X^\delta_s+z)}\psi^{(1)}_s\chi'_s}
  {|k|^{-1/2}+i\frac{k}{|k|}\cdot V^\delta_s}\\
&
  \frac{(X_s-X_s^\delta)\cdot(V_s-V_s^\delta)}{|\delta|^{4/3}|X_s-X_s^\delta|}
  \:dk\:dz\:ds,
\end{split}
\]
correspondingly
\[\begin{split}
& \mbox{IV}(t,x,v)=  \int_0^t\int_{\RR^d}\int_{\RR^d}\tilde{\alpha}(k)\:\frac{k\cdot
    z}{|k|\:|z|^{d}}\:\frac{e^{ik\cdot(X^\delta_s+z)}\chi_s}
  {|k|^{-1/2}+i\frac{k}{|k|}\cdot V^\delta_s}\\
  &\left[-\nabla_3\psi^{(1)}_s
    \frac{|z|}{|X_s-X_s^\delta|^3}(X_s-X^\delta_s)\cdot(V_s-V^\delta_s)\right.\\
  &  \left.
    +\nabla_2\psi^{(1)}_s\cdot\left(\frac{V_s-V^\delta_s}{|X_s-X^\delta_s|}
      -\frac{X_s-X^\delta_s}{|X_s-X^\delta_s|^3}
      (X_s-X_s^\delta)\cdot(V_s-V_s^\delta)\right)\right]
  \:dk\:dz\:ds,
\end{split}\]
and
\[
\mbox{V}(t,x,v)= \int_0^t\int_{\RR^{2d}}\tilde{\alpha}(k)\:\frac{k\cdot
    z}{|k|\:|z|^{d}}\:\frac{e^{ik\cdot(X^\delta_s+z)}\chi_s\psi^{(1)}_s}
  {\left(|k|^{-1/2}+i\frac{k}{|k|}\cdot V^\delta_s\right)^{2}}\:
  i\frac{k}{|k|}\cdot F(X^\delta_s) \:dk\:dz\:ds.
\] 
Let us define
\begin{equation*}
%   \label{eq:def-I}
  \mbox{I}(T)=\int_0^T\iint_{\Omega}
  \frac{V_t-V^\delta_t}{A_\delta(t,x,v)}
  \cdot \mbox{I}(t,x,v) \:dx\:dv\:dt,
\end{equation*}
and $\mbox{II}(T)$, $\mbox{III}(T)$, $\mbox{IV}(T)$ and $\mbox{V}(T)$
similarly.

We are going to bound each of these terms. The last one gives the good
order of $\int_{\RR^d}|k|^{3/2+a}|\alpha(k)|^2\:dk$. The others
are bounded by integrals involving lower powers of $|k|$.

\subsubsection{Upper bound for $|\textmd{V}(T)|$}
\label{sec:V}

First, we make the change of variables $z'=z+X^\delta_s$, followed by
the change of variable $(x',v')=(X_s^\delta,V_s^\delta)$ in the
integral defining $\mbox{V}(T)$. When $(x,v)\in\Omega$, the variable
$(x',v')$ belongs to the set
$\Omega_s=\{(X^\delta(s,x,v),V^\delta(s,x,v)),\:(x,v)\in\Omega\}$.
Note also that $X(-s,x',v')=x+\delta_1$ and
$V(-s,x',v')=v+\delta_2$.

Writing for convenience $x,v,z$ instead of $x',v',z'$, it follows
from these changes of variables and from Property~\ref{prop:Hamilt}
that
\begin{multline*}
  \mbox{V}(T)=\int_0^T\int_0^t\iint_{\Omega_s}
  \int_{\RR^d}\int_{\RR^d}
  \tilde{\chi}_s\tilde{\psi}^{(1)}_s \\
  \begin{aligned}
    \frac{\tilde{V}^\delta_{t,s}-V_{t-s}}{|\delta|^2
      +\sup_{0\leq r\leq t}|\tilde{X}^\delta_{r,s}-X_{r-s}|^2
      +\int_0^t|\tilde{V}^\delta_{r,s}-V_{r-s}|^2\:dr}\cdot \tilde{\alpha}(k) \\
    \frac{k\cdot(z-x)}{|k|\:|z-x|^d}\:e^{ik\cdot z}
    \frac{i\frac{k}{|k|}\cdot F(x)}
    {\left(|k|^{-1/2}+i\frac{k}{|k|}\cdot v\right)^2}\:dk\:dz\:dx\:dv\:ds\:dt,
  \end{aligned}
\end{multline*}
where
\begin{gather*}
  \tilde{X}^\delta_{t,s}=X^\delta(t,X(-s,x,v),V(-s,x,v)), \\
  \tilde{V}^\delta_{t,s}=V^\delta(t,X(-s,x,v),V(-s,x,v)), \\
  \tilde{\psi}^{(1)}_s=\psi^{(1)}\left(\frac{z-x}{|z-x|},
    \frac{\tilde{X}^\delta_{s,s}-x}{|\tilde{X}^\delta_{s,s}-x|},
    \frac{|z-x|}{|\tilde{X}^\delta_{s,s}-x|}\right).
\end{gather*}
and
\begin{equation*}
  \tilde{\chi}_s=\chi\left(
\frac{|\tilde{X}^\delta_{s,s}-x|}{|\delta|^{4/3}}\right).
\end{equation*}
Writing the tensor (remember that $\alpha(k)\in\RR^d$)
\begin{equation*}
  G_{\text{V}}(v,z)=\int_{\RR^d}\frac{k\otimes
      k}{|k|^2}\otimes\frac{\tilde{\alpha}(k) \: e^{ik\cdot
        z}}{\left(|k|^{-1/2}+i\frac{k}{|k|}\cdot v\right)^2}\:dk,
\end{equation*}
and reminding that $\Omega_s\subset\Omega'$ for all $s\in[0,T]$, we
have
\begin{align*}
  |\text{V}(T)| & \leq C\int_0^T\iint_{\Omega'}
  \int_0^t\int_{\RR^d}
  \frac{\tilde{\chi}_s\tilde{\psi}^{(1)}_s
    \:|F(x)|\:\|G_{\text{V}}(v,z)\|}{|z-x|^{d-1}
    \left(|\delta|+|\tilde{X}^\delta_{s,s}-x|\right)} \\ & \qquad\qquad
  \frac{|\tilde{V}^\delta_{t,s}-V_{t-s}|}{\left(|\delta|^2+
      \int_0^t|\tilde{V}^\delta_{r,s}-V_{r-s}|^2\:dr\right)^{1/2}}
  \:dz\:ds\:dx\:dv\:dt,
\end{align*}
where $\|a\|^2=\sum_{i,j,k=1}^d a_{ijk}^2$ for any tensor with three
entries $a=(a_{ijk})$ with ${1\leq i,j,k\leq d}$. So
\begin{align*}
  |\text{V}(T)|
 & \leq C
  \int_0^T\iint_{\Omega'}\int_{\RR^d}
  \frac{\tilde{\psi}^{(1)}_s\:\|G_{\text{V}}(v,z)\|}{|z-x|^{d-1}
    \left(|\delta|+|\tilde{X}^\delta_{s,s}-x|\right)} \\ & \qquad\qquad
  \int_s^T\frac{|\tilde{V}^\delta_{t,s}-V_{t-s}|}{\left(|\delta|^2+
      \int_0^t|\tilde{V}^\delta_{r,s}-V_{r-s}|^2\:dr\right)^{1/2}}\:dt
  \:dz\:dx\:dv\:ds.
\end{align*}
Now, on the one hand, following the same computation as in
Lemma~\ref{lem:1}, the integral with respect to $t$ can be upper
bounded by $C(\log 1/|\delta|)^{1/2}$ for $|\delta|$ small enough.  On
the other hand,
\begin{align*}
  & \iint_{\Omega'}\int_{\RR^d}
  \frac{\tilde{\psi}^{(1)}_s\:\|G_{\text{V}}(v,z)\|}{|z-x|^{d-1}
    \left(|\delta|+|\tilde{X}^\delta_{s,s}-x|\right)} \:dz\:dx\:dv
  \notag \\ & \qquad \leq
  \left(\iint_{\Omega'}\int_{\RR^d}
    \frac{\tilde{\psi}^{(1)}_s}{|z-x|^{d-1}
      \left(|\delta|+|\tilde{X}^\delta_{s,s}-x|\right)}
    \:dz\:dx\:dv\right)^{1/2} \notag \\ & \qquad\qquad\qquad
  \left(\int_{\Omega'_2}\int_{\RR^d}\int_{\Omega_1}
    \frac{\tilde{\psi}^{(1)}_s\:\|G_{\text{V}}(v,z)\|^2}{|z-x|^{d-1}
      \left(|\delta|+|\tilde{X}^\delta_{s,s}-x|\right)}
    \:dx\:dz\:dv\right)^{1/2}, \notag
\end{align*}
and this last term is bounded by
\begin{align}
 & C
  \left(\iint_{\Omega'}\frac{1}{|\tilde{X}^\delta_{s,s}-x|}
    \int_{x+K(\tilde{X}^\delta_{s,s}-x)}
    \frac{dz}{|z-x|^{d-1}}
    \:dx\:dv\right)^{1/2} \notag \\ & \qquad\quad
  \left(\int_{\Omega'_2}\int_{\RR^d}\|G_{\text{V}}(v,z)\|^2\int_{\Omega'_1}
    \frac{dx}{|z-x|^{d-1}
      \left(|\delta|+|z-x|\right)}
    \:dz\:dv\right)^{1/2} \notag \\ & \qquad\qquad \qquad\quad
  \leq C(\log 1/|\delta|)^{1/2}
  \left(\int_{\Omega'_2}\int_{\RR^d}\|G_{\text{V}}(v,z)\|^2
    \:dz\:dv\right)^{1/2}, \label{eq:bourrin}
\end{align}
where we have used that, for any $z\in x+K(\tilde{X}^\delta_{s,s}-x)$,
$|z-x|\leq|\tilde{X}^\delta_{s,s}-x|$, and where the last inequality
can be obtained by a spherical coordinate change of variable centered
at $x$ in the variable $z$ in the first term, and centered at $z$ in
the variable $x$ in the second term.

Now,
\begin{align*}
  & \int_{\Omega'_2}\int_{\RR^d}\|G_{\text{V}}(v,z)\|^2 \:dz\:dv
  \notag \\
  & \qquad =\int_{\Omega'_2}\int_{\RR^d}\sum_{i,j,n=1}^d\int_{\RR^d}\int_{\RR^d}
  \frac{k_il_ik_jl_j}{|k|^2|l|^2}
  \frac{\tilde{\alpha}_n(k)\overline{\tilde{\alpha}_n(l)}e^{iz\cdot(k-l)}}
  {\left(|k|^{-1/2}+i\frac{k}{|k|}\cdot
      v\right)^2}\\
&\qquad\qquad\qquad\left(|l|^{-1/2}-i\frac{l}{|l|}\cdot
      v\right)^{-2}\:dl\:dk\:dz\:dv \notag,
\end{align*}
and integrating first in $z$ and $l$, this is equal to
\[ 
  \int_{\Omega'_2}\sum_{i,j}^d\int_{\RR^d}
  \frac{k_i^2k_j^2}{|k|^4}
  \frac{|\tilde{\alpha}(k)|^2}{\left||k|^{-1/2}+i\frac{k}{|k|}\cdot
      v\right|^4}\:dk\:dv. \notag
\]
Therefore
\begin{align*}
 \int_{\Omega'_2}\int_{\RR^d}\|G_{\text{V}}(v,z)\|^2 \:dz\:dv & 
  \leq C\int_{\RR^d}\int_{\Omega'_2}
  \frac{|\tilde{\alpha}(k)|^2}{\left(\frac{1}{|k|}+\left(\frac{k\cdot
          v}{|k|}\right)^2\right)^2}\:dv\:dk \notag \\ & 
  \leq C\int_{\RR^d}|k|^2|\tilde{\alpha}(k)|^2\int_{-\infty}^{+\infty}
  \frac{dv_1}{(1+|k|v_1^2)^2}\:dk, \notag
\end{align*}
where we  write the vector $v$ as
$(v_1,\ldots,v_d)$ in an orthonormal basis of $\RR^d$ with first
vector $k/|k|$. In conclusion
\begin{align}
 \int_{\Omega'_2}\int_{\RR^d}\|G_{\text{V}}(v,z)\|^2 \:dz\:dv
 & 
  \leq C\int_{\{|k|>(\log 1/|\delta|)^2\}}|k|^{3/2}|\alpha(k)|^2\:dk
  \notag \\ & \leq C(\log
  1/|\delta|)^{-4a}\int_{\RR^d}|k|^{\frac{3}{2}+2a}|{\alpha}(k)|^2\:dk,
  \label{eq:Plancherel} 
\end{align}
Combining this inequality with~(\ref{eq:bourrin}), we finally get
\begin{equation*}
  |\text{V}(T)|\leq C
  (\log 1/|\delta|)^{1-2a}\left(\int_{\RR^d}
    |k|^{\frac{3}{2}+2a}|{\alpha}(k)|^2\:dk\right)^{1/2}.
\end{equation*}

\subsubsection{Upper bound for $|\textmd{IV}(T)|$}
\label{sec:IV}

Applying to $\text{IV}(T)$ the same change of variable as we did for
$\text{V}(T)$, we have
\begin{multline*}
  |\text{IV}(T)|\leq C\int_0^T\iint_{\Omega'}
  \int_0^t\int_{\RR^d}
  \frac{|\tilde{V}^\delta_{s,s}-v|\:|\tilde{V}^\delta_{t,s}-V_{t-s}|}
  {|\delta|^2+\int_0^t|\tilde{V}^\delta_{r,s}-V_{r-s}|^2\:dr}
  \\ \frac{\tilde{\chi}_s\:|\hat{\psi}^{(1)}_s|\:\|G_{\text{IV}}(v,z)\|}
  {|z-x|^{d-1}|\tilde{X}^\delta_{s,s}-x|}\:dz\:ds\:dx\:dv\:dt
\end{multline*}
where $\|a\|^2=\sum_{i,j=1}^d a_{ij}^2$ for any matrix $a=(a_{ij})_{1\leq i,j\leq d}$,
\begin{equation*}
  G_{\text{IV}}(v,z)=\int_{\RR^d}\frac{k}{|k|}
\otimes\:\frac{\tilde{\alpha}(k) \: e^{ik\cdot
      z}}{|k|^{-1/2}+i\frac{k}{|k|}\cdot v}\:dk
\end{equation*}
and
\begin{multline*}
  \hat{\psi}^{(1)}_s=-\nabla_3\psi^{(1)}\left(\frac{z-x}{|z-x|},
    \frac{\tilde{X}^\delta_{s,s}-x}{|\tilde{X}^\delta_{s,s}-x|},
    \frac{|z-x|}{|\tilde{X}^\delta_{s,s}-x|}\right)
  \frac{|z|\:(\tilde{X}^{\delta}_{s,s}-x)}{|\tilde{X}^{\delta}_{s,s}-x|^2}
  \\
  \begin{aligned}
    -\left(\nabla_2\psi^{(1)}\left(\frac{z-x}{|z-x|},
        \frac{\tilde{X}^\delta_{s,s}-x}{|\tilde{X}^\delta_{s,s}-x|},
        \frac{|z-x|}{|\tilde{X}^\delta_{s,s}-x|}\right)\cdot
      \frac{\tilde{X}^{\delta}_{s,s}-x}{|\tilde{X}^{\delta}_{s,s}-x|}\right)
    \frac{\tilde{X}^{\delta}_{s,s}-x}{|\tilde{X}^{\delta}_{s,s}-x|} \\
    +\nabla_2\psi^{(1)}\left(\frac{z-x}{|z-x|},
      \frac{\tilde{X}^\delta_{s,s}-x}{|\tilde{X}^\delta_{s,s}-x|},
      \frac{|z-x|}{|\tilde{X}^\delta_{s,s}-x|}\right).
  \end{aligned}
\end{multline*}
Note that, because of the properties of $\psi^{(1)}$ obtained in
Section~\ref{sec:CV-z}, 
\begin{equation*}
  |\hat{\psi}^{(1)}_s|\leq C \mathbb{I}_{\{z-x\in K(\tilde{X}^{\delta}_{s,s}-x)\}}
\end{equation*}
for some constant $C$.

Then, following a similar computation as the one leading to~(\ref{eq:bourrin}),
\begin{align*}
  |\text{IV}(T)| & \leq C
  \left(\int_{\Omega'_2}\int_{\RR^d}\|G_{\text{IV}}(v,z)\|^2\int_{\Omega'_1}
    \frac{dx\:dz\:dv}{|z-x|^{d-1}
      \left(|\delta|^{4/3}+|z-x|\right)}
    \right)^{1/2} \notag \\
  & \qquad\qquad
  \left(\iint_{\Omega'}\int_0^T\int_0^t
    \frac{|\tilde{V}^\delta_{s,s}-v|^2\:|\tilde{V}^\delta_{t,s}-V_{t-s}|^2}
    {|\delta|^4+\left(\int_0^t|\tilde{V}^\delta_{r,s}-V_{r-s}|^2\:dr\right)^2}
  \right. \notag \\  & \qquad\qquad\qquad
  \left.\frac{1}{|\tilde{X}^\delta_{s,s}-x|}
    \int_{x+K(\tilde{X}^\delta_{s,s}-x)}\frac{dz}{|z-x|^{d-1}}
    \:ds\:dt\:dx\:dv\right)^{1/2} \notag. 
\end{align*}
Hence
\begin{multline}
  |\text{IV}(T)| \leq C(\log 1/|\delta|)^{1/2}
  \left(\int_{\Omega'_2}\int_{\RR^d}\|G_{\text{IV}}(v,z)\|^2\right)^{1/2}
  \\
  \left(\iint_{\Omega'}\int_0^T\int_0^t
    \frac{|\tilde{V}^\delta_{s,s}-v|^2\:|\tilde{V}^\delta_{t,s}-V_{t-s}|^2}
    {|\delta|^4+\left(\int_0^t|\tilde{V}^\delta_{r,s}-V_{r-s}|^2\:dr\right)^2}
  \right)^{1/2} \label{eq:bourrin-2}
\end{multline}
where we have used that $|\tilde{X}^\delta_{s,s}-x|\geq
|\delta|^{4/3}$ and $|\tilde{X}^\delta_{s,s}-x|\geq|z-x|$ when
$\tilde{\chi}_s\:|\hat{\psi}^{(1)}_s|\not =0$.

Now, making the change of variable
$(x',v')=(X^\delta(-s,x.v),V^\delta(-s,x,v))$ and denoting $(x',v')$
as $(x,v)$ for convenience, we have
\begin{align*}
  \iint_{\Omega'}\int_0^T\int_0^t & 
  \frac{|\tilde{V}^\delta_{s,s}-v|^2\:|\tilde{V}^\delta_{t,s}-V_{t-s}|^2}
  {|\delta|^4+\left(\int_0^t|\tilde{V}^\delta_{r,s}-V_{r-s}|^2\:dr\right)^2}
  \notag \\
  & \leq \iint_{\Omega''}\int_0^T
  \frac{|V_t-V^\delta_t|^2\int_0^t|V_s-V^\delta_s|^2\:ds}
  {|\delta|^4+\left(\int_0^t|V_s-V^\delta_s|^2\:ds\right)^2}\:dt\:dx\:dv
  \notag \\
  & =\frac{1}{2}\iint_{\Omega''}\log\left(
    \frac{|\delta|^4+\left(\int_0^T|V_s-V^\delta_s|^2\:ds\right)^2}{|\delta|^4}
  \right)\:dx\:dv \notag \\ & \leq C\log(1/|\delta|).
%   \label{eq:V_s-V_t}
\end{align*}
Next, similarly as in the computation leading to~(\ref{eq:Plancherel}), we have
\begin{align*}
  \int_{\Omega'_2}\int_{\RR^d}\|G_{\text{IV}}(v,z)\|^2 \:dz\:dv 
  & \leq C\int_{\RR^d}\int_{\Omega'_2}
  \frac{|\tilde{\alpha}(k)|^2}{\frac{1}{|k|}+\left(\frac{k\cdot
          v}{|k|}\right)^2}\:dv\:dk \\ &
  \leq C\int_{\RR^d}|k|\,|\tilde{\alpha}(k)|^2\int_{-\infty}^{+\infty}
  \frac{dv_1}{1+|k|v_1^2}\:dk \\ & 
  \leq C\int_{\RR^d}|k|^{1/2}|\tilde{\alpha}(k)|^2\:dk
  \\ & \leq C(\log
  1/|\delta|)^{-2-4a}\int_{\RR^d}|k|^{\frac{3}{2}+2a}|{\alpha}(k)|^2\:dk.
\end{align*}

The combination of these inequalities finally yields
\begin{equation*}
  |\text{IV}(T)| \leq C\left(\int_{\RR^d}
    |k|^{\frac{3}{2}+2a}|{\alpha}(k)|^2\:dk\right)^{1/2}
\end{equation*}
if $|\delta|<1/e$.

\subsubsection{Upper bound for $|\textmd{III}(T)|$}
\label{sec:III}

As before, we compute
\begin{multline*}
  |\text{III}(T)|\leq C\int_0^T\iint_{\Omega'}
  \int_0^t\int_{\RR^d}
  \frac{|\tilde{V}^\delta_{s,s}-v|\:|\tilde{V}^\delta_{t,s}-V_{t-s}|}
  {|\delta|^2+\int_0^t|\tilde{V}^\delta_{r,s}-V_{r-s}|^2\:dr}
  \\ \frac{\tilde{\psi}^{(1)}_s\:|\tilde{\chi}'_s|\:\|G_{\text{IV}}(v,z)\|}
  {|\delta|^{4/3}|z-x|^{d-1}|}\:dz\:ds\:dx\:dv\:dt.
\end{multline*}

Then, proceeding as in~(\ref{eq:bourrin-2}),
\begin{align*}
  |\text{III}(T)| & \leq \frac{C}{|\delta|^{4/3}}
  \left(\int_{\Omega'_2}\int_{\RR^d}\|G_{\text{IV}}(v,z)\|^2\int_{B(z,2|\delta|^{4/3})}
    \frac{dx}{|z-x|^{d-1}}
    \:dz\:dv\right)^{1/2} \\
  & \qquad\qquad
  \left(\iint_{\Omega'}\int_0^T\int_0^t
    \frac{|\tilde{V}^\delta_{s,s}-v|^2\:|\tilde{V}^\delta_{t,s}-V_{t-s}|^2}
    {|\delta|^4+\left(\int_0^t|\tilde{V}^\delta_{r,s}-V_{r-s}|^2\:dr\right)^2}
  \right. \\  & \qquad\qquad
  \left.\mathbb{I}_{\{|\tilde{X}^\delta_{s,s}-x|\leq 2|\delta|^{4/3}\}}
    \int_{x+K(\tilde{X}^\delta_{s,s}-x)}\frac{dz}{|z-x|^{d-1}}
    \:ds\:dt\:dx\:dv\right)^{1/2}, 
\end{align*}
so that
\begin{align*}
  |\text{III}(T)|
 & \leq C(\log 1/|\delta|)^{1/2}
  \left(\int_{\Omega'_2}\int_{\RR^d}\|G_{\text{IV}}(v,z)\|^2\right)^{1/2}
\end{align*}
where we have used that $|z-x|\leq|\tilde{X}^\delta_{s,s}-x|\leq
2|\delta|^{4/3}$ when $\tilde{\psi}^{(1)}_s\:|\tilde{\chi}'_s|\not=0$.

Finally,
\begin{equation*}
  |\text{III}(T)| \leq C\left(\int_{\RR^d}
    |k|^{\frac{3}{2}+2a}|{\alpha}(k)|^2\:dk\right)^{1/2}
\end{equation*}
if $|\delta|<1/e$.

\subsubsection{Upper bound for $|\textmd{I}(T)|$ and $|\textmd{II}(T)|$}
\label{sec:I-II}

We only detail the computation of a bound for
$|\text{I}(T)|$. The case of $|\text{II}(T)|$ is very similar and is
left to the reader.

We compute as before
\begin{multline*}
  |\text{I}(T)|\leq C\int_0^T\iint_{\Omega'}
  \int_{\RR^d}
  \frac{|\tilde{V}^\delta_{t,t}-v|}
  {\left(|\delta|^2+\int_0^t|\tilde{V}^\delta_{r,t}-V_{r-t}|^2\:dr\right)^{1/2}}
  \\ \frac{\tilde{\chi}_s\:\tilde{\psi}^{(1)}_s\:\|G_{\text{IV}}(v,z)\|}
  {|z-x|^{d-1}\big(|\delta|+|\tilde{X}^\delta_{t,t}-x|\big)}\:dz\:dx\:dv\:dt.
\end{multline*}

Next, the computation is very similar to~(\ref{eq:bourrin-2}):
\begin{align*}
  |\text{I}(T)| & \leq C
  \left(\int_{\Omega'_2}\int_{\RR^d}\|G_{\text{IV}}(v,z)\|^2\int_{\Omega'_1}
    \frac{dx}{|z-x|^{d-1}(|\delta|+|z-x|)}
    \:dz\:dv\right)^{1/2} \\
  & \qquad\qquad\qquad
  \left(\iint_{\Omega'}\int_0^T
    \frac{|\tilde{V}^\delta_{t,t}-v|^2}
    {|\delta|^2+\int_0^t|\tilde{V}^\delta_{r,t}-V_{r-t}|^2\:dr}
  \right. \\  & \qquad\qquad\qquad\qquad
  \left.\frac{1}{|\tilde{X}^\delta_{t,t}-x|}
    \int_{x+K(\tilde{X}^\delta_{t,t}-x)}\frac{dz}{|z-x|^{d-1}}
    \:dt\:dx\:dv\right)^{1/2}.
\end{align*}
Proceeding as before
\begin{align*}
  |\text{I}(T)|
& \leq C(\log 1/|\delta|)^{1/2}
  \left(\int_{\Omega'_2}\int_{\RR^d}\|G_{\text{IV}}(v,z)\|^2\right)^{1/2} \\
  & \qquad\qquad\qquad
  \left(\iint_{\Omega''}\int_0^T
    \frac{|V_t-V^\delta_t|^2\:dt}
    {|\delta|^2+\int_0^t|V_r-V^\delta_{r}|^2\:dr}
  \:dx\:dv\right)^{1/2},
\end{align*}
so that eventually
\begin{align*}
  |\text{I}(T)|
 & \leq C\log (1/|\delta|)
  \left(\int_{\Omega'_2}\int_{\RR^d}\|G_{\text{IV}}(v,z)\|^2\right)^{1/2}
  \\ & \leq C\left(\int_{\RR^d}
    |k|^{\frac{3}{2}+2a}|{\alpha}(k)|^2\:dk\right)^{1/2}.
\end{align*}
\medskip

This completes the proof that 
\begin{equation*}
  |B^{11}_\delta(T)|\leq C\left(1+(\log
    |1/|\delta|)^{1-2a}\right)\left(\int_{\RR^d} 
    |k|^{\frac{3}{2}+2a}|{\alpha}(k)|^2\:dk\right)^{1/2},
\end{equation*}
where
\begin{equation*}
  B^{11}_\delta(T):=\int_0^T\iint_{\Omega}
  \frac{V_t-V_t^\delta}{A_\delta(t,x,v)}\cdot\int_0^t\chi_s
  \:B^{11}_\delta(s,x,v)\:ds\:dx\:dv\:dt.
\end{equation*}
%%%%%%%%%%%%%%%%%%%%%%%%%%%%%%%%
\subsubsection{Upper bound for $|B^{12}_\delta(T)|$}
\label{sec:B^12}
%%%%%%%%%%%%%%%%%%%%%%%%%%%%%%%%%
Let us define
\begin{equation*}
  B^{12}_\delta(T):=\int_0^T\iint_{\Omega}
  \frac{V_t-V_t^\delta}{A_\delta(t,x,v)}\cdot\int_0^t\chi_s
  \:B^{12}_\delta(s,x,v)\:ds\:dx\:dv\:dt.
\end{equation*}
As will appear below, this term is very similar to
$\text{V}(T)$. 

We apply the same method as before, {\em without integrating by parts in
time}:
\begin{multline*}
  |B^{12}_\delta(T)|\leq C\int_0^T\iint_{\Omega'}
  \int_0^t\int_{\RR^d}
  \frac{|\tilde{V}^\delta_{t,s}-V_{t-s}|}
  {\left(|\delta|^2+\int_0^t|\tilde{V}^\delta_{r,s}-V_{r-s}|^2\:dr\right)^{1/2}}
  \\ \frac{\tilde{\psi}^{(1)}_s\:\|G_{12}(v,z)\|}
  {|z-x|^{d-1}(|\delta|+|\tilde{X}^\delta_{s,s}-x|)}\:dz\:ds\:dx\:dv\:dt
\end{multline*}
where
\begin{equation*}
  G_{12}(v,z)=\int_{\RR^d}\frac{k}{|k|^{1/2}}
\otimes\:\frac{\tilde{\alpha}(k) \: e^{ik\cdot
      z}}{|k|^{-1/2}+i\frac{k}{|k|}\cdot v}\:dk.
\end{equation*}
Again,
\begin{align*}
  |B^{12}_\delta(T)| & \leq C\log(1/|\delta|)
  \left(\int_{\Omega'_2}\int_{\RR^d}\|G_{12}(v,z)\|^2\:dz\:dv\right)^{1/2}
  \\ & \leq C\log(1/|\delta|)
  \left(\int_{\Omega'_2}\int_{\RR^d}\frac{|k|\,|\tilde{\alpha}(k)|^2}
    {\frac{1}{|k|}+\left(\frac{k\cdot
          v}{|k|}\right)^2}\:dk\:dv\right)^{1/2}
  \\ & \leq C\log(1/|\delta|)
  \left(\int_{\RR^d}|k|^{3/2}|\tilde{\alpha}(k)|^2\:dk\right)^{1/2}
  \\ & \leq C(\log 1/|\delta|)^{1-2a}
  \left(\int_{\RR^d}|k|^{\frac{3}{2}+2a}|\tilde{\alpha}(k)|^2\:dk\right)^{1/2}.
\end{align*}

\subsubsection{Conclusion}
\label{sec:ccl}

Combining all the previous inequalities, we obtain that
\begin{equation*}
  |B^{1}_\delta(T)|\leq C\left(1+(\log 1/|\delta|)^{1-2a}\right)
  \left(\int_{\RR^d}|k|^{\frac{3}{2}+2a}|\tilde{\alpha}(k)|^2\:dk\right)^{1/2}
\end{equation*}
where
\begin{equation*}
  B^{1}_\delta(T):=\int_0^T\iint_{\Omega}
  \frac{V_t-V_t^\delta}{A_\delta(t,x,v)}\cdot\int_0^t\chi_s
  \:B^{1}_\delta(s,x,v)\:ds\:dx\:dv\:dt.
\end{equation*}

Now, we observe from~(\ref{eq:def-B1-B2}) that $B^{2}(s,x,v)$ has
exactly the same structure as $B^{1}(s,x,v)$: a singularity of order
$d-1$ in $z$, a function $\psi^{(2)}$ that has all the required
regularity, and a term $e^{ik\cdot(X^\delta_s+z)}$. It is then easy to
see that this term can be treated by exactly the same method as
$B^{1}(s,x,v)$. We leave the details to the reader.

From this follows that
\begin{equation*}
  |B_\delta(T)|\leq C\left(1+(\log 1/|\delta|)^{1-2a}\right)
  \left(\int_{\RR^d}|k|^{\frac{3}{2}+2a}|\tilde{\alpha}(k)|^2\:dk\right)^{1/2}
\end{equation*}
where
\begin{equation*}
  B_\delta(T):=\int_0^T\iint_{\Omega}
  \frac{V_t-V_t^\delta}{A_\delta(t,x,v)}\cdot\int_0^t\chi_s
  \:B_\delta(s,x,v)\:ds\:dx\:dv\:dt.
\end{equation*}

Finally, the term $C_\delta(s,x,v)$ of~(\ref{eq:def-C}) can be bounded
exactly as $B_\delta(s,x,v)$ by simply exchanging the roles of $X_s$
and $X_s^\delta$. Therefore, the proof of
Proposition~\ref{prop:control-Q} is completed.

%
%%%%%%%%%%%%%%%%%%%%%%%%%%%%%%%%%%%%%%%%%%%%%%%%%%

\end{document}